%% file: main.tex
\documentclass{scrartcl}

\input {definitions.tex}

\begin {document}

\title {Signed counts of\\ real simple rational functions}
\author {Boulos El Hilany, Johannes Rau
{\footnote{{For this work, the first and second author were supported 
by the DFG Research Grant RA 2638/2-1.}}}
{}
}

\maketitle

\input {Introduction.tex}
\input {dessins.tex}
\input {bipartitegraphs.tex}

\input {invariance.tex}
\input {vanishing.tex}

\input {brokensequences.tex}

\input {polynomiality.tex}
\input {nonvanishing.tex}

\printbibliography

\subsection*{Contact}

\begin{itemize}
	\item 
  Boulos El Hilany,
  Instytut Matematyczny Polskiej Akademii Nauk,
  ul. \'Sniadeckich 8,
	00-656 Warszawa,
	Poland;
\href{boulos.hilani@gmail.com}{boulos.hilani@gmail.com}.

  \item
  Johannes Rau, 
	Universität Tübingen,
  Geschwister-Scholl-Platz, 
	72074 Tübingen, 
	Germany;
  \href{mailto:johannes.rau@math.uni-tuebingen.de}{johannes.rau@math.uni-tuebingen.de}.
\end{itemize}

\subsection*{Acknowledgement}

We would like to thank Frédéric Bihan, Ilia Itenberg, Hannah Markwig and Arthur Renaudineau for numerous helpful discussions.
We would also like two thank the two anonymous referees for careful reading and many useful suggestions and corrections.

\end {document}

%% file: definitions.tex
\KOMAoptions{paper=a4}
\KOMAoptions{DIV=10} 

\usepackage[utf8]{inputenc} 
\usepackage[T1]{fontenc} 
\usepackage[ngerman,english]{babel} 
\usepackage{csquotes} 
\usepackage{todonotes} 

\usepackage{lmodern}
\usepackage{avant}
\usepackage{courier}

\usepackage{enumitem}
\setlist[enumerate]{label*=(\alph*),ref=(\alph*)}

\KOMAoptions{DIV=last}

\usepackage {amsmath} 
\usepackage{amssymb} 
\usepackage{amscd} 

\usepackage{bm}
\usepackage{float}
\usepackage{multirow}

\newcommand{\even}{{\text{even}}}
\newcommand{\odd}{{\text{odd}}}
\DeclareMathOperator{\Sym}{Sym}
\DeclareMathOperator{\dis}{dis}
\DeclareMathOperator{\lev}{lev}
\DeclareMathOperator{\Rot}{Rot}
\DeclareMathOperator{\Tr}{Tr}
\DeclareMathOperator{\Flip}{Flip}
\DeclareMathOperator{\sech}{sech}
\DeclareMathOperator{\Res}{Res}
\DeclareMathOperator{\pol}{pol}

\DeclareMathOperator{\Des}{Des}
\DeclareMathOperator{\conj}{conj}

\newcommand{\R}{\mathbb{R}}

\newcommand{\Fatlambda}{\mathbf{\Lambda}}
\newcommand{\fatlambda}{{\bm \lambda}}

\newcommand{\cc}{\mathfrak{c}}
\newcommand{\oo}{\mathfrak{o}}
\newcommand{\ee}{\mathfrak{e}}
\newcommand{\bb}{\mathfrak{b}}
\newcommand{\SBB}{\mathbf{SB}_{(B)}}
\newcommand{\SBC}{\mathbf{SB}_{(C)}}
\newcommand{\BWG}{E}

\renewcommand{\sp}{\text{s}}

\newcommand{\CC}{{\mathbf C}}

\newcommand{\NN}{{\mathbf N}}
\newcommand{\PP}{{\mathbf P}}
\newcommand{\QQ}{{\mathbf Q}}
\newcommand{\RR}{{\mathbf R}}

\newcommand{\ZZ}{{\mathbf Z}}

\newcommand{\MMM}{{\mathcal M}}
\newcommand{\PPP}{{\mathcal P}}

\newcommand{\SSS}{{\mathcal S}}

\usepackage{graphicx} 
\graphicspath{{pic/}} 

\usepackage[format=plain,labelfont={sf,footnotesize},textfont=small,margin=12pt]{caption} 
\DeclareCaptionLabelSeparator{MySpace}{\enskip } 
\usepackage{booktabs} 

\usepackage[style=alphabetic,
						sorting=nyt,
						maxnames=4,
						backend=biber,
						isbn=false,
						clearlang=true]
						{biblatex}
\bibliography{biblio}
\DeclareFieldFormat*{title}{\textit{#1}} 
\renewbibmacro{in:}{} 
\DeclareFieldFormat{journaltitle}{\iffieldequalstr{journaltitle}{ArXiv e-prints}{Preprint}{#1}} 
\DeclareFieldFormat{booktitle}{#1} 
\DeclareRedundantLanguages{english,English}{english,german,ngerman,french}
\DeclareSourcemap{ 
  \maps[datatype=bibtex]{
    \map{
      \step[fieldsource=volume,final]
      \step[fieldset=doi,null]
			\step[fieldset=url,null]
    }  
  }
}

\usepackage{hyperref} 
\usepackage{xcolor} 
\definecolor{darkblue}{RGB}{0,0,170}
\definecolor{darkred}{RGB}{200,0,0}
\hypersetup{citecolor=darkblue, urlcolor=darkblue, linkcolor=darkblue, colorlinks=true} 
\usepackage[all]{hypcap} 

\usepackage{amsthm} 

\newtheorem {theorem}{Theorem}[section]
\newtheorem {proposition}[theorem]{Proposition}
\newtheorem {lemma}[theorem]{Lemma}

\newtheorem {corollary}[theorem]{Corollary}

\theoremstyle {definition}
\newtheorem {definition}[theorem]{Definition}
\newtheorem {example}[theorem]{Example}

\newtheorem {construction}[theorem]{Construction}

\theoremstyle {remark}
\newtheorem {remark}[theorem]{Remark}
\newtheorem {convention}[theorem]{Convention}

\newtheorem {notation}[theorem]{Notation}

%% file: Introduction.tex
\begin{abstract}
	We study the problem of counting real simple rational functions $\varphi$
	with prescribed ramification data (i.e.\ a particular class
	of oriented real Hurwitz numbers of genus $0$). 
	We introduce a signed count of such functions 
	which is independent of the position of the branch points,
	thus providing a lower bound for the actual count (which does depend on the position). 
	We prove (non-)vanishing theorems
	for these signed counts and study
	their asymptotic growth when adding further simple branch points.
	The approach is based on \cite{IZ18} which 
	treats the polynomial case.
\end{abstract}

\section{Introduction}

A \emph{simple rational function} of degree $d$ is a function $\varphi : \CC\PP^1 \to \CC\PP^1$ 
which, in affine coordinates, has the form
\[
  \varphi(z) = \frac{\psi(z)}{z-p}
\]
with $\psi \in \CC[z], \deg(\psi) = d, p \in \CC$ and $\psi(p) \neq 0$.
We call $\varphi$ \emph{real} if $\psi \in \RR[z]$ and $p \in \RR$,
and \emph{increasing} if the leading coefficient of $\psi$ is positive.
Two such functions are considered equivalent if they differ 
by linear coordinate change $z \mapsto \lambda z + \mu, \lambda, \mu \in \RR, \lambda > 0$.
A unique representative in each equivalence class of increasing functions
is given by \emph{normalized} functions, 
which are of the form $\varphi(z) = \psi(z)/z$, where $\psi(0) \neq 0$ and 
the leading coefficient of $\psi$ is $1$.

Let $b_1, \dots, b_k \in \CC$ be the critical values (also called branch points)
of $\varphi$ and let $\Lambda_1, \dots, \Lambda_k$ be the corresponding ramification profiles. 
This means $\Lambda_j = (\Lambda_j^1 \geq \dots \geq \Lambda_j^{l_j})$ is the partition of $d$
such that $b_j$ has $l_j$ preimages at which $f$ locally takes the form
$z \mapsto z^{\Lambda_j^i}$. 
We call $\{(b_j, \Lambda_j)\}_{j=1,\dots,k}$ the \emph{ramification data} of $\varphi$.

Let us now fix $k$ distinct \emph{real} points $\PPP = (b_1, \dots, b_k), b_j \in \RR$ and 
$k$ partitions of $d$ $\Fatlambda= (\Lambda_1, \dots, \Lambda_k)$ and set 
$d = \sum_{i,j} (\Lambda_j^i - 1)$.
We are interested in the \emph{number of normalized real simple rational
functions $f$ of degree $d$ with ramification data $\{(b_j, \Lambda_j)\}_{j=1,\dots,k}$}.
We denote the set of such functions by $\SSS(\PPP, \Fatlambda)$.
The number $|\SSS(\PPP, \Fatlambda)|$ does not depend on the position of the branch points, as long
as we do not change their order of appearance on the real line.
However, in general $|\SSS(\PPP, \Fatlambda)|$ is \emph{not} invariant under a permutation of 
this order.
To remedy the situation, we will define a sign $\varepsilon(\varphi) \in \{\pm 1\}$ for any $\varphi\in \SSS(\PPP, \Fatlambda)$ (see Definition \ref{def_signsimplerationalfunction})
such that the following theorem holds.

\begin{theorem}[Invariance theorem] \label{Th:InvarTh}
  The number 
	\[
	  \sum_{\varphi\in \SSS(\PPP, \Fatlambda)} \varepsilon(\varphi)
	\]
	does not depend on the position of the branch points $b_1, \dots, b_k \in \RR$.
	We call it the \emph{$S$-number} of $\Fatlambda$ and denote it by 
	$S(\Fatlambda)$.
\end{theorem}

Note that by definition $|S(\Fatlambda)|$ gives a lower bound for $|\SSS(\PPP, \Fatlambda)|$.
It is therefore interesting to find criteria for $S(\Fatlambda) \neq 0$ (equivalent to the existence of functions)
or to prove statements about the asymptotic growth.
To formulate our results, we need the following notation. 
For a partition $\Lambda = (\Lambda^1, \dots, \Lambda^t)$ of $d$, 
the \emph{reduced} partition $\lambda$ is the partition obtained from $(\Lambda^1 - 1, \dots, \Lambda^t - 1)$
after removing all zeros. A branch point is called \emph{simple} if its reduced ramification profile is $\lambda = (1)$. 
Let us now fix $k$ decreasing finite sequences $\fatlambda=(\lambda_1, \dots, \lambda_k)$ and $m \in \NN$. We set 
\[
  d := m + \sum_{i,j} \lambda_j^i, \vspace{-1.5ex}
\]
and let $\Fatlambda = (\Lambda_1, \dots, \Lambda_{k+m})$ be the unique collection of partitions of $d$ 
such that $\Lambda_1, \dots, \Lambda_k$ reduce to $\lambda_1, \dots, \lambda_k$ and 
$\Lambda_{k+1} = \dots = \Lambda_{k+m}$ reduce to $(1)$. 
Hence $S(\fatlambda, m) := S(\Fatlambda)$ is the $S$-number of normalized simple rational functions with
$k+m$ branch points, $k$ of which have reduced ramification profiles $\lambda_1, \dots, \lambda_k$
while the remaining $m$ branch points are simple. 
We collect these numbers in two generating series, separately for $d$ odd and even:
\[
  F^\odd_\fatlambda(q) = \sum_{\substack{m = 0 \\ d \text{ odd}}}^\infty S(\fatlambda, m) \frac{q^m}{m!}, \quad \quad
  F^\even_\fatlambda(q) = \sum_{\substack{m = 0 \\ d \text{ even}}}^\infty S(\fatlambda, m) \frac{q^m}{m!}.
\]

\begin{theorem}[(Non-)Vanishing theorem] \label{Th:IffOdd} \label{Th:IffEven} \label{thm_nonvanishing}
  The generating series $F^\odd_\fatlambda(q)$ is not identically zero if and only if the following conditions hold:
	\begin{itemize}
		\item \vspace{-0.8ex} In each partition $\lambda_j$ at most one odd number appears an odd number of times and at most one even number appears an odd number of times.
		\item \vspace{-0.8ex} There exists an even number of partitions $\lambda_j$ having exactly one even number appearing an odd number of times.
	\end{itemize}
	\vspace{-0.8ex} The generating series $F^\even_\fatlambda(q)$ is not identically zero if and only if the following conditions hold:
	\begin{itemize}
		\item \vspace{-0.8ex} In each partition $\lambda_j$ at most one odd number appears an odd number of times and at most one even number appears an odd number of times.
	\end{itemize}
\end{theorem}

\begin{theorem}[Logarithmic growth] \label{thm_LogGrowth}
  Fix $\fatlambda=(\lambda_1, \dots, \lambda_k)$ and the parity of $d$ such that the non-vanishing criteria
	in Theorem~\ref{thm_nonvanishing} are satisfied. Set $\pi \equiv d - \sum_{i,j} \lambda_j^i \mod 2$ the corresponding
	parity of $m$. Then the logarithmic growth of $S(\fatlambda, m)$ in $m$ is given by
	\[
	  \ln |S(\fatlambda, m)| \underset{\substack{m \to \infty \\ m \equiv \pi \hspace{-1.5ex} \mod 2}}{\sim} m \ln(m).
	\]
\end{theorem}

These statements should be compared to \cite[Theorems 1,3,4,5]{IZ18} in the polynomial case.

Let $H^\CC(\fatlambda, m)$ denote the Hurwitz number counting
complex simple rational functions $\varphi \in \CC(z)$
with $k$ critical levels of reduced ramification type $\lambda_1, \dots, \lambda_k$ and 
$m$ additional simple branch points.
Note that $\ln H^\CC(\fatlambda, m) \sim m \ln(m)$ for $m \to \infty$ (see Remark \ref{rem_complexcount}), so 
under the conditions of the theorem the real and complex counts are logarithmically equivalent.

Since all the main theorems crucially depend on the definition of the signs $\varepsilon(\varphi)$,
let us insert its definition here.

\begin{definition} 
  Let $\Sigma = (a_1, \dots, a_n)$ be a finite sequence of integers. 
	A \emph{disorder} of $\Sigma$ is a pair $i < j$ such that $a_i > a_j$. 
	The number of disorders is denoted by $\dis(\Sigma)$. 
\end{definition}

For any $a \in \CC$ we denote by $R_\varphi(a) \in \NN$ the ramification index
of $\varphi$ at $a$, i.e.\ the order of vanishing of $\varphi(z) - \varphi(a)$ for $a \neq p$, and $R_\varphi(p) = 1$.

\begin{definition} \label{def_signsimplerationalfunction}
  Let $\varphi$ be a real simple rational function with simple pole $p$ and let $b$ be a branch point of $\varphi$. 
  We set $\Sigma_b$ to be the sequence of ramification indices $R_\varphi(a)$ for 
	$a \in (\varphi^{-1}(b) \cap \RR) \sqcup \{p\}$, ordered according to their appearance on the real line.
	Let $b_1, \dots, b_k$ denote the collection of branch points of $\varphi$.
	The \emph{sign} of $\varphi$ is
	\[
	  \varepsilon(\varphi) := (-1)^{\dis(\varphi)} \quad \quad 
		\text{where } \dis(\varphi) := \dis(\Sigma_{b_1}) + \dots + \dis(\Sigma_{b_k}).
	\]
\end{definition}

\begin{remark}
	The main difference with \cite{IZ18} is that our definition considers the simple pole $p$ 
	as part of the preimage for any critical level. 
	It will become clear in Section \ref{sec_bipartitegraphs} that this is, under some assumptions,
	the only definition with a chance to satisfy to Theorem \ref{Th:InvarTh}.
	We call the disorders involving $p$ \emph{pole disorders} 
	and all other ones \emph{level disorders}.
\end{remark}

Let us give a brief outline of the paper. 
In Sections \ref{sec_dessins}, \ref{sec_bipartitegraphs} and \ref{sec_invariance}
we closely follow the approach in \cite{IZ18} to prove the Invariance Theorem 
\ref{Th:InvarTh} using dessins d'enfant. The main difference in these sections
is that instead of trees we deal with graphs with a loop.
Our presentation focuses on the difference to \textit{loc.~cit.}~while 
being mostly self-contained.
Sections \ref{sec vanishing}, \ref{Sec:Pol} and \ref{Sec:Non-Van} deal with Theorems \ref{thm_nonvanishing}
and \ref{thm_LogGrowth}.
Here, the differences to \textit{loc.~cit.}~are more significant. 
In particular, in Section~\ref{Sec:BrokenAlt} we introduce the generating series of 
\emph{broken alternations}, show that it obeys a certain differential equation
and use this to express it in terms of the generating series
of (ordinary) alternations. 
This is crucial for the extension of the ideas from 
\textit{loc.~cit.}~to the case of simple rational functions.

The counts of simple rational functions under investigation here 
can be considered as \emph{oriented} versions of 
\emph{real Hurwitz numbers} as defined for example in \cite{MR3, GMR}. 
The approach to study these numbers via an invariant signed count
is in the spirit of Welschinger invariants \cite{We, IKS4}, even though
the definition of Welschinger signs is of different flavour than 
Definition \ref{def_signsimplerationalfunction}. 
Real Hurwitz numbers and their invariance/asymptotic properties
have also been studied in the context of topological quantum field theories, 
matrix models, moduli spaces of real algebraic curves \cite{AN, GZ3, Orl}, 
however, mostly in the context of \emph{completely imaginary} 
configurations of branch points $\PPP$, in contrast 
to the \emph{completely real} configurations $\PPP \subset \RR\PP^1$ we consider here.

%% file: dessins.tex
\section{Dessins d'enfant for simple rational functions} \label{sec_dessins}

Conjugation-invariant dessins d'enfant are used to  describe
real functions in 
\cite{Bar-SpaceRealPolynomials, NSV-TopologicalClassificationGeneric, IZ18}
(for generic polynomials, generic rational functions and arbtirary polynomials, respectively).
The following definition is adapted to describe the simple rational functions 
considered in this paper.

We fix $d$ and a sequence $\Fatlambda = (\Lambda_1,\ldots,\Lambda_k)$ of partitions of $d$ such that
$d = \sum_{i,j} (\Lambda_j^i - 1)$. Note that this implies $k > 1$.
For a partition $\Lambda$ we denote by $l(\Lambda)$ the number of parts of $\Lambda$.
For the vertex $v$ of a graph $G$, its \emph{degree} $\deg(v)$ is the number of adjacent half edges.
Throughout the following, we fix an affine chart $\CC\PP^1 = \CC \sqcup \{\infty\}$ 
and use the orientation on $\RR\PP^1$ induced by the increasing orientation on $\RR$.

\begin{definition}\label{Def:RealSimpleDessin}
A \emph{real simple rational dessin} of degree $d$ and type $\Fatlambda$ is a graph $\Gamma \subset \CC\PP^1$ 
whose vertices are labelled by elements of the set $\{1,\ldots,k,\infty\}$ such that the following conditions hold.
\begin{enumerate}
	\item The labelled graph $\Gamma$ is invariant under complex conjugation $\conj$.
	\item The real circle $\RR\PP^1$ is a union of edges of $\Gamma$.
	\item Exactly two vertices of $\Gamma$ are labelled by $\infty$, namely $P:= \infty \in \CC\PP^1$ and $p \in \RR$, 
	      with $\deg(P) = 2d-2$ and $\deg(p)=2$.
	\item For each integer $1\leq j\leq k$, the graph $\Gamma$ has exactly $l(\Lambda_j)$ 
	      vertices labelled by $j$ and their degrees are equal to the elements of $\Lambda_j$, multiplied by $2$.
	\item Each edge of $\Gamma$ is one of the following $k+1$ types: $\infty\rightarrow 1$, $1\rightarrow 2$, $2\rightarrow 3$, $\ldots$, $k-1\rightarrow k$, $k\rightarrow\infty$;
	      in particular, this induces an orientation on $\Gamma$.
	\item For any connected component $C$ of $\CC\PP^1 \setminus \Gamma$, each type of edges appears exactly once in the boundary $\partial C$ of $C$.
	\item There exists an edge in $\RR\PP^1$ of type $k \to \infty$ whose orientation agrees with that of $\RR\PP^1$.
\end{enumerate}
\end{definition}

\begin{remark}
  Given a real dessin $\Gamma \subset \CC\PP^1$, we call $\Gamma \setminus \{P\} \subset \CC$ the \emph{affine dessin} for $\Gamma$. 
	It is a graph with $2d - 2$ unbounded ends of type $\infty \to 1$ or $k \to \infty$, and it is easy to adapt the above
	conditions into \emph{affine} versions such that both contain exactly the same information. We will mostly use affine dessins in our figures.
	An example for $d=21, k=4$ is given in Figure \ref{Fig:GraphDessin}. The vertex $p$ is drawn in red.
\end{remark}

It is easy to associate a real simple rational dessin $\Gamma_\varphi$ to any increasing real simple rational function $\varphi$
with branch points $b_1 < \dots < b_k \in \RR$.
As a set, we define $\Gamma_\varphi := \varphi^{-1}(\RR\PP^1)$. We label the two preimages of $\infty \in \CC\PP^1_{\text{target}}$ by $\infty$ 
and the preimages of $b_j$ by $j$. 

\begin{definition}
  Two real simple rational dessins $\Gamma, \Gamma'$ are \emph{equivalent} if there exists 
  a homeomorphism $\Phi : \CC\PP^1 \to \CC\PP^1$ such that 
	$\Phi$ and $\conj$ commute, 
	$\Phi|_{\RR\PP^1}$ is orientation-preserving,
	$\Phi(\Gamma) = \Gamma'$ and $\Phi$ preserves the labels.
	We denote the set of equivalence classes by $\Des(\Fatlambda)$.
\end{definition}

\begin{theorem} \label{thm_rationalfunctionasdessins}
  Fix arbitrary real points $\PPP = \{b_1 < \dots < b_k\} \in \RR\PP^1$. Then the map
  \[
	  \SSS(\PPP, \Fatlambda) \to \Des(\Fatlambda), \quad \varphi \mapsto [\Gamma_\varphi],
	\]
	is a bijection.
\end{theorem}

\begin{proof}
  The proof is a straightforward application of dessins d'enfant and the Riemann existence theorem
	and can be easily adapted from \cite[Propositions 2.6 and 2.7]{IZ18}.
\end{proof}

\begin{remark} 
  We define the sign of dessin to be the sign 
	of the associated rational function, $\epsilon(\Gamma_\varphi) = \epsilon(\varphi)$.
	Note that $\epsilon(\Gamma_\varphi)$ can be easily read of from the dessin itself
	since the sequence $\Sigma_{b_j}$ of ramification indices for the branch point $b_j$ 
	is equal to the sequence of degrees of the vertices labelled
	by $i$ plus the pole $p$ divided by $2$ (cf.\ Definition \ref{def_signsimplerationalfunction}).
\end{remark}

%% file: bipartitegraphs.tex
\section{Bipartite graphs} \label{sec_bipartitegraphs}

In this section, we introduce auxiliary combinatorial objects that are used in the proof of the invariance theorem. 
They are equivalent to simple rational dessins with $k=2$ critical values.

\begin{definition}\label{Def:BWGraphs} \label{Def:DisorderSignsGraphs}
	A \emph{black and white simple graph}, or short, \emph{$bw$-graph} is a 
	connected graph $G$ embedded into $\CC$ whose vertices are 
	coloured in black and white in alternation, and it has only one cycle of edge length $2$ 
	(see Figure~\ref{Fig:TwoGraphExamples}). 
	A $bw$-graph is said to be \emph{real} if it is invariant (including the colours) 
	under complex conjugation. 
	Two (real) $bw$-graphs are isomorphic if one can be transformed into the other by 
	an ($\conj$-equivariant) homeomorphism of $\CC$ (that preserves the orientation
	of $\RR$).

	For a real $bw$-graph $G$ the \emph{real part} is $\RR G := G \cap \RR$. 
	The leftmost and rightmost vertices of $\RR G$ are called the \emph{border vertices}.
	The graph is called \emph{white-sided} (respectively \emph{black-sided}) if its rightmost 
	border vertex is white (respectively, black).
	It is called \emph{short} if the the cycle of $G$ contains a border vertex.
	Otherwise, we call it \emph{long}.

\begin{figure}[b]
\centering
\includegraphics[scale=1.4]{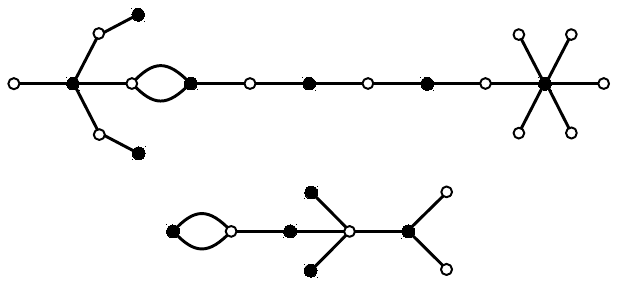}

\caption{Two real $bw$-graphs.}

\label{Fig:TwoGraphExamples}
\end{figure}

	We denote by $(\Sigma_w,\Sigma_b)$ the sequence of degrees of its real white and real black vertices, 
	respectively, from left to right.
	A \emph{level disorder} is a disorder in one of these sequences, while a pole disorder 
	is real vertex of $G$ to the left of the cycle
	and with degree larger than $1$ (including the left cycle vertex). 
 	In this section, it is more convenient to keep track of
	level and pole disorders separately. 
	We denote them by $\lev(G)$ and $\pol(G)$, respectively.
	The \emph{sign} of $G$ is 
	\[
	  \varepsilon (G)=(-1)^{\lev(G) + \pol(G)}.
	\]
\end{definition}

\begin{example}\label{Ex:RealVert}\label{Ex:SignGraph}
The real $bw$-graph appearing at the top of  Figure~\ref{Fig:TwoGraphExamples} is a white-sided long graph. 
It has $(\Sigma_w,\Sigma_b)=(132221,43226)$ and its sign is equal to 
$(-1)^{\displaystyle 12 + 2}=+1$. \\
The real $bw$-graph $G$ appearing at the bottom of Figure~\ref{Fig:TwoGraphExamples} 
is a black-sided short graph.
It has $(\Sigma_w,\Sigma_b)=(34,223)$ and its sign is equal to 
$(-1)^{\displaystyle 0 + 1}=-1$.
\end{example}

\begin{remark}\label{Rem:012or3}
  A vertex of $G$ is called  \emph{even} or \emph{odd} if its degree $\deg(v)$ is even or odd, respectively. 
	A real $bw$-graph $G$ can have four, two or zero real odd vertices.
	Only the first case 
  is equivalent to $G$ being long. Note that in this case not all four vertices have the same colour. 
\end{remark}

For the rest of this section, we fix a positive integer $d$, and two partitions $\Lambda_w,\Lambda_b$ of $d$. 
We denote by $\bm{G}$ the set of real $bw$-graphs such that the degrees of their white and black vertices give $\Lambda_w$ and $\Lambda_b$, respectively. 
We split $\bm{G} = \bm{W} \sqcup \bm{B}$ into the subsets of white-sided and black-sided graphs, respectively.

\begin{convention}
  Let $M$ be a set together with a sign function $\varepsilon : M \to \{\pm 1\}$.
	We set
	\[
		S(M) := \sum_{x \in M} \varepsilon(x).
	\]
\end{convention}

\begin{theorem}\label{Th:InvThGraphs} Fix a positive integer $d$ and two partitions $\Lambda_w$ and $\Lambda_b$ of $d$. Then $$ S(\bm{W}) =S(\bm{B}).$$ 
\end{theorem} 

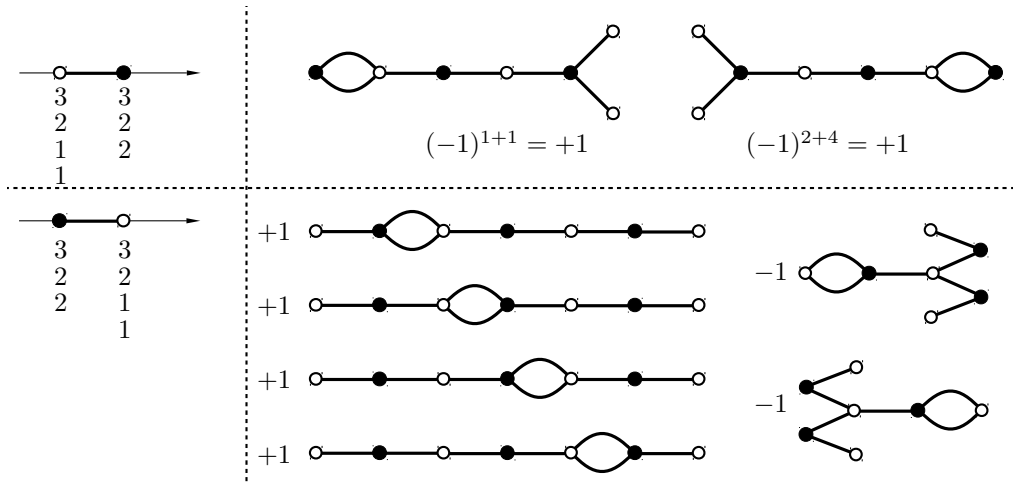
\begin{figure}
\centering
\includegraphics[scale=1.4]{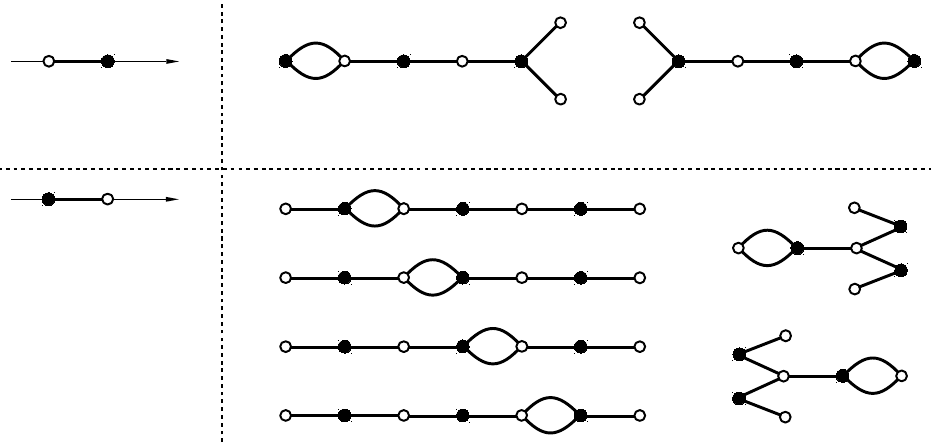}
\put(-335,143){\small{$3$}}
\put(-335,133){\small{$2$}}
\put(-335,123){\small{$2$}}
\put(-359,143){\small{$3$}}
\put(-359,133){\small{$2$}}
\put(-359,123){\small{$1$}}
\put(-359,113){\small{$1$}}
\put(-220,125){\small{$(-1)^{1+1}=+1$}}
\put(-100,125){\small{$(-1)^{2+4}=+1$}}
\put(-359,85){\small{$3$}}
\put(-359,75){\small{$2$}}
\put(-359,65){\small{$2$}}
\put(-335,85){\small{$3$}}
\put(-335,75){\small{$2$}}
\put(-335,65){\small{$1$}}
\put(-335,55){\small{$1$}}
\put(-283,92){\small{$+1$}}
\put(-283,64){\small{$+1$}}
\put(-283,36){\small{$+1$}}
\put(-283,8){\small{$+1$}}
\put(-97,27){\small{$-1$}}
\put(-97,77){\small{$-1$}}
\caption{An example of $S(\bm{W}) =S(\bm{B})$ for $\Lambda_w=(3,2,1,1)$ and $\Lambda_b=(3,2,2)$.}
\label{Fig:ExampleInvar}
\end{figure}

\begin{example}\label{Ex:InvTh}
Figure~\ref{Fig:ExampleInvar} shows all graphs with $\Lambda_w=(3,2,1,1)$ and $\Lambda_b=(3,2,2)$. 
The sum of signs is equal to 2 both for the black-sided graphs on top and for the white-sided graphs on the bottom.
\end{example}

\subsection{Symmetrizing graphs} 

To prove Theorem~\ref{Th:InvThGraphs}, we subsequently identify subsets of graphs which cancel each other,
hence reducing the problem to special classes of graphs.
The first step is to symmetrize graphs in a certain sense.

\begin{figure}[h]
\centering
\includegraphics[scale=1.3]{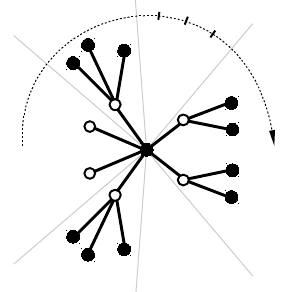}
\put(-90,50){\small{$T_1$}}
\put(-82,105){\small{$T_2$}}
\put(5,50){\small{$T_k$}}
\caption{Ordering of the pairs of non-real trees growing at $v$.}
\label{Fig:TreeOrder}
\end{figure}

\begin{definition} \label{def_pairoftrees}
	Let $v$ be a real vertex of the graph $G$.  A \emph{pair of trees $T$ growing at $v$} is a pair of 
	complex conjugated connected components $T = (\Theta, \overline{\Theta})$ of $G \setminus \{v\}$ 
	with $\Theta \neq \overline{\Theta}$. 
	Let $T_1, \dots, T_k$ be the sequence of all pairs of trees growing at $v$,
	ordered according to their appearance in $\{z : \Im(z) > 0\}\subset\CC$
	following the clock-wise direction
	(see Figure \ref{Fig:TreeOrder}). 
	Then the \emph{forest} growing at $v$ is the sequence	$F(v) = (T_1, \dots, T_{k-1})$ 
	if $v$ is an odd border vertex and $F(v) = (T_1, \dots, T_{k})$ otherwise.
\end{definition}

Let $l(v)$ denote the number of entries of $F(v)$. 
If $v$ is even, then $\deg(v) = 2k + 2 = 2 l(v) + 2$. If $v$ is odd, then $\deg(v) = 2k + 1$ or $\deg(v) = 2k + 3$ depending on whether $v$ is a border vertex or not. Due to the shortening of $F(v)$ for border vertices we obtain $\deg(v) = 2 l(v) + 3$ for any odd vertex with $\deg(v) \neq 1$.

\begin{figure}
    \centering
    \begin{minipage}{.5\textwidth}
        \centering
        \includegraphics[scale=1]{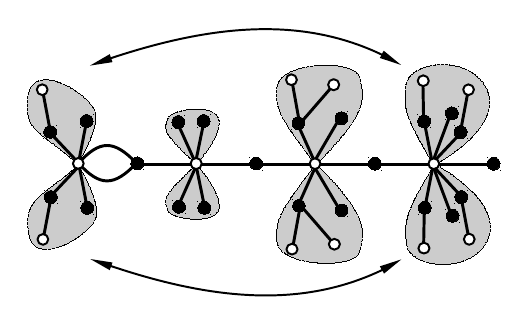}
    \end{minipage}
    \begin{minipage}{0.49\textwidth}
        \centering
        \includegraphics[scale=1.8]{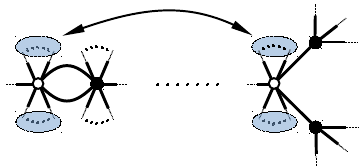}
    \end{minipage}
  	\caption{Exchanging forests for even and odd vertices, respectively.}
		\label{Fig:ExchangForests}
\end{figure}

\begin{definition}\label{Def:flip}
  A collection of real vertices $v_1, \dots, v_k$ of a given graph $G$ are \emph{of the same type} if they are of the same colour, the same parity and $\deg(v_i) > 1$ for all $i$.
Given a permutation $\sigma$ of the $v_i$, we obtain a new graph $\sigma(G)$ by cutting off the forest $F(v_i)$ from $v_i$ and replanting it at $v_{\sigma(v_i)}$, for all $i$. Hereby we keep the order of the individual pairs of trees in $F(v_i)$ and the position of the (fixed) right-most pair of trees in the case of a odd border vertex (see Figure~\ref{Fig:ExchangForests}). This determines $\sigma(G)$ uniquely up to homeomorphism. In the case of two vertices $v_1, v_2$ and $\sigma$ the transposition, we call the result the \emph{flipped graph} $\Flip_{v_1}^{v_2}(G)$.
\end{definition}

Note that in terms of the real sequences of $G$ the operation $\Flip_{v_1}^{v_2}$ acts like a transposition permuting the two entries corresponding to $v_1$ and $v_2$.

\begin{definition}\label{Def:EvSymm}
Let $\Sigma=(e_1, \dots, e_{2m})$ be an integer sequence $\Sigma$ of even length $2m$. Then $\Sigma$ is \emph{symmetric} if $e_{m-k} = e_{m+k+1}$ for all $k = 0,\dots, m-1$.
If $\Sigma$ is not symmetric, the \emph{first non-symmetric pair} denotes the unique pair $(m-k,m+k+1)$ such that $e_{m-k} \neq e_{m+k+1}$ and 
$(e_{m-k+1}, \dots, e_{m+k})$ is symmetric.
 
 A (general) integer sequence $\Sigma$ is \emph{nearly symmetric} if it is symmetric after applying the following two steps:
\begin{itemize}
	\item Remove all $1$'s.
	\item If the resulting sequence is of odd length, remove the first entry.
\end{itemize}
If $\Sigma$ is not nearly symmetric then the \emph{first non-symmetric pair} is the pair of indices which corresponds to the first non-symmetric
pair after applying the same two steps.

Given an integer sequence $\Sigma$, we denote by $\Sigma^{\even}$ and $\Sigma^{\odd}$ the subsequence of even and odd entries, respectively. 
A real $bw$-graph $G$ is \emph{nearly symmetric} if all the sequences $\Sigma^{\even}_w, \Sigma^{\odd}_w, \Sigma^{\even}_b, \Sigma^{\odd}_b$ are nearly symmetric separately.  
If $G$ is not nearly symmetric, then the \emph{first non-symmetric pair} is the pair of vertices corresponding to the first non-symmetric pair
for the first sequence from before (in the given order) which is not nearly symmetric.
\end{definition}

\begin{remark} \label{rem_nearlysymmetric}
  Let us unravel the previous definition in the given situation. Let $\Sigma^{\even}$ be a nearly symmetric sequence of only even entries, then  either it is of the form $(e_1, \dots, e_m, e_m, \dots, e_1)$ or $(e_0, e_1, \dots, e_m, e_m, \dots, e_1)$. Let $\Sigma^{\odd}$ be a sequence of only odd entries and of length at most $3$	and not of the form $(1,1)$ or $(x,1,y)$ (which cannot  appear in practice). Then it is nearly symmetric if and only if it is of one of the following forms (where we assume $a \neq 1 \neq b$).
	\[
		\emptyset, (a), (1), (a,a), (a,1), (1,a), (1,a,1), (a,a,1), (1,a,a), (a,b,b)
	\]
	If $\Sigma^{\odd}$ is not nearly symmetric, it is of the form $(\underline{c},\underline{d})$, $(\underline{c},\underline{d},1)$, 
	$(1,\underline{c},\underline{d})$, $(a,\underline{c},\underline{d})$ with $c\neq d$. The first non-symmetric pair is formed by the underlined entries.
\end{remark}

\begin{example} 
  The top graph in Figure~\ref{Fig:TwoGraphExamples} is not nearly symmetric and its first non-symmetric pair
	is formed by the left- and right-most black vertex of degree $4$ and $6$, respectively.
	The bottom graph is nearly symmetric. 
\end{example}

\begin{lemma}\label{lem_flipsignchange}
  Let $G$ be a graph which is not nearly symmetric and let $v_1 < v_2$ be the first non-symmetric pair of $G$. Then
	\[
	  \varepsilon(\Flip_{v_1}^{v_2} (G)) = - \varepsilon(G).
	\]
\end{lemma}

\begin{proof}
Let $\Sigma$ be the sequence of degrees of some colour and parity as $v_1,v_2$. If $\deg v_1$ and $\deg v_2$ are consecutive entries of $\Sigma$,
flipping the two obviously changes the number of disorders by $1$. If not consecutive, they are separated by a stretch of the form
$(e_1, \dots, e_m, e_m, \dots, e_1)$. Such a stretch can be removed without changing the parity of $\dis(\Sigma)$,
since any disorder involving the stretch shows up an even number of times. Hence the previous case applies. 
The pole disorders are not affected by the flipping.
\end{proof}

\begin{proposition} \label{pro_symmetrizinggraphs}
  Consider the subsets $\bm{N\!W}\subset\bm{W}$ and $\bm{N\! B}\subset\bm{B}$ of nearly symmetric white-sided and black-sided graphs, respectively. 
	Then we have
	\begin{equation}\label{eq:Sum=Sum}
		S(\bm{W}) = S(\bm{N\!W}), \quad \quad S(\bm{B}) = S(\bm{N\!B}).
	\end{equation} 
\end{proposition}

\begin{proof}
  We define the involution $\Sym : \bm{W} \setminus \bm{N\!W} \to \bm{W} \setminus \bm{N\!W}$ 
	by setting $\Sym(G) = \Flip_{v_1}^{v_2}(G)$, where $v_1 < v_2$ is the first non-symmetric pair of $G$. 
	By Lemma~\ref{lem_flipsignchange} the pairs $G, \Sym(G)$ cancel out in $S(\bm{W})$
	and the statement follows.
\end{proof}

\subsection{Rotating and shifting vertices} 

The next step is based on rotating graphs by $180^\circ$. This operation reverses the real sequences.
In order to stay in the class of nearly symmetric graphs, we need to make a slight adjustment.

\begin{figure}[H]
\centering
\includegraphics[scale=0.9]{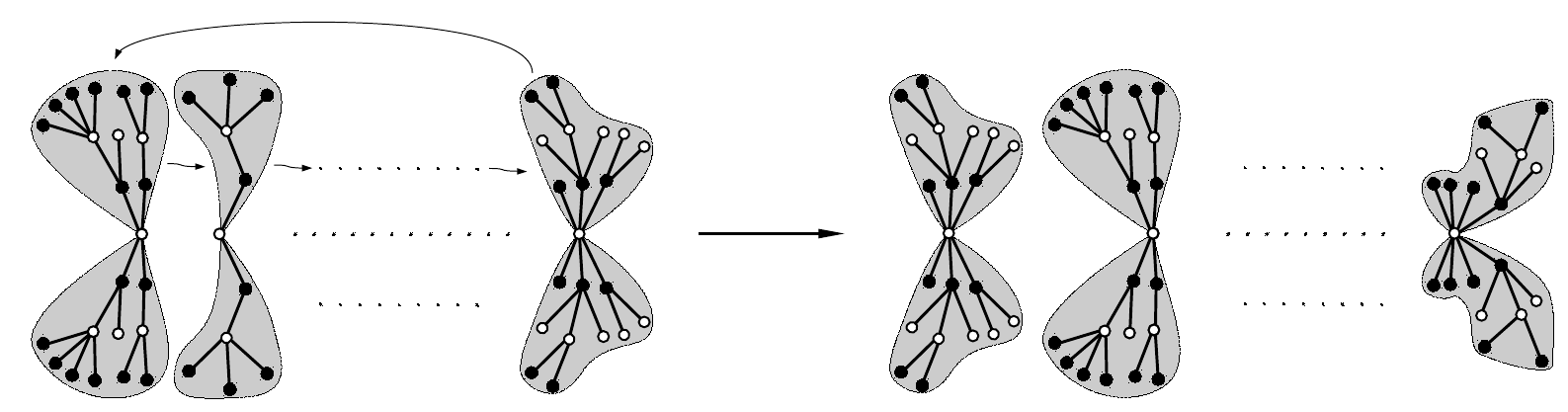}
\put(-400,-10){$F_1$}
\put(-360,-10){$F_2$}
\put(-270,-10){$F_k$}
\put(-180,-10){$F_k$}
\put(-120,-10){$F_1$}
\put(-15,-10){$F_{k-1}$}
\caption{The cyclic shift operation.}
\label{Fig:CycledShifting}
\end{figure}
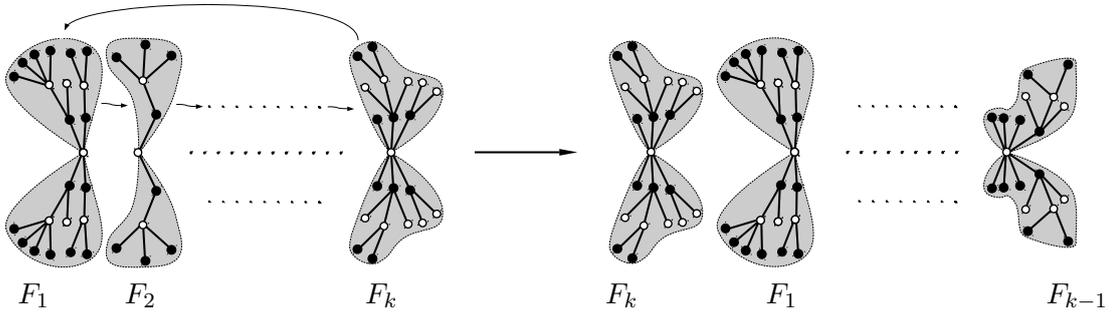

\begin{definition}\label{Def:CycledShif}
Let $v_1 < \ldots < v_k$ be a collection of ordered real vertices of the same type and let $\sigma = (v_1 \dots v_k)$ be the cyclic permutation of these vertices. 
The \emph{cyclic shift} of $G$ based on $v_1 < \ldots < v_k$ is the graph $\sigma(G)$ (see Definition~\ref{Def:flip} and Figure~\ref{Fig:CycledShifting}).
\end{definition}

\begin{remark}\label{Rem:Restoring}
	Let $\Sigma$ be a nearly symmetric sequence and let $\Sigma'$ be the reversed sequence. 
	Then either $\Sigma'$ is nearly symmetric or can be made nearly symmetric by 
	cyclically shifting all non-one entries of $\Sigma'$ to the next non-one entry (keeping all $1$'s fixed). 
	Obviously, if $\Sigma'$ is part of the real sequence of a graph $G$, then the shifting, if necessary, can be performed by the cycled shifting operation
	based on the vertices corresponding to the non-one entries of $\Sigma'$. 
\end{remark} 

We set $\bm{N\!G} = \bm{N\!W} \cup \bm{N\!B}$.

\begin{definition}\label{Def:Rotation}
We define the map $ \Rot:~\bm{N\!G} \to \bm{N\!G}$, where $\Rot(G)$ is the graph obtained from $G$ by a rotation of $180^\circ$ and, if necessary, performing cyclic shifts based  
on the higher-valent vertices of certain colour and parity
(see Remark~\ref{Rem:Restoring}). 
\end{definition}

\begin{lemma} \label{lem_rotationsignchange}
  Let $G$ be a nearly symmetric graph and let $\Sigma$ be its real sequence of white (or black) vertices. 
	Let $\Rot(\Sigma)$ denote the real sequence of white (or black) vertices of $\Rot(G)$. Set $l := l(\Sigma^{\even})$.	
	Then we have
	\[
	  \dis(\Rot(\Sigma)) - \dis(\Sigma) \equiv
		\begin{cases}
		  0 & \Sigma^{\odd} = \emptyset \text{ or }(a,a), \\ 
		  1 & \Sigma^{\odd} = (1,a) \text{ or }(a,1), \\
			l  & \text{otherwise},
		\end{cases} \mod 2.
	\]
	Here we refer to the ten possible forms for $\Sigma^{\odd}$ introduced in Remark \ref{rem_nearlysymmetric}.
\end{lemma}

\begin{proof}
  We set $\overline{\Sigma} = \Rot(\Sigma)$.
	As in Lemma \ref{lem_flipsignchange}, throughout the proof we will use the fact that a consecutive pair
	of identical entries in $\Sigma$ can be removed without changing the parity of $\dis(\Sigma)$. 
	We first focus on $\Sigma^{\odd}$. Going through the ten types in Remark \ref{rem_nearlysymmetric} we can
	easily check that $\dis(\overline{\Sigma}^{\odd}) - \dis(\Sigma^{\odd}) \equiv 1$ only if $\Sigma^{\odd} = (1,a)$ or $(a,1)$
	and $\dis(\overline{\Sigma}^{\odd}) - \dis(\Sigma^{\odd}) \equiv 0$ otherwise. 
	
	We now reduce to this case by the following trick: Let $k$ be the number of transposition necessary to 
	transform $\Sigma$ into $(\Sigma^{\odd}, \Sigma^{\even})$ (i.e.\ moving all odd entries to the beginning of the sequences).
	Since each of these transpositions flips an odd and even entry (i.e.,\ two non-equal numbers), we get
	$\dis(\Sigma) \equiv \dis((\Sigma^{\odd}, \Sigma^{\even})) + k$. The pattern of odd and even entries in $\overline{\Sigma}$
	is reversed to the pattern of $\Sigma$ (the cyclic shifts only permute odd and even entries separately). 
	Hence moving all odd entries in $\overline{\Sigma}$ to the right, we get
	$\dis(\overline{\Sigma}) \equiv \dis((\overline{\Sigma}^{\even}, \overline{\Sigma}^{\odd})) + k$. 
	If we set $m := l(\Sigma^{\odd})$ and move all the odd entries back to the left, we get 
	$\dis(\overline{\Sigma}) \equiv \dis((\overline{\Sigma}^{\odd}, \overline{\Sigma}^{\even})) + k + lm$.
	Putting things together, we obtain
	\begin{align} 
	  \nonumber
	  \dis(\Rot(\Sigma)) - \dis(\Sigma) &\equiv \dis((\overline{\Sigma}^{\odd}, \overline{\Sigma}^{\even})) - \dis((\Sigma^{\odd}, \Sigma^{\even})) + lm \\
		                                  &\equiv \dis(\overline{\Sigma}^{\odd}) - \dis(\Sigma^{\odd}) + lm.
	\end{align}
	Plugging in the possible types for $\Sigma^{\odd}$, the claim follows.
\end{proof}

\begin{lemma}\label{L:RotSafeGraph}
Let $G\in\bm{N\!G}$ be a nearly symmetric graph. We have $\varepsilon(\Rot(G))=-\varepsilon(G)$ if there exists a sequence $\Sigma_c^{\odd}(G)$ of $G$ of some colour $c$ such that $\Sigma_c^{\odd}(G)=(1,b,b)$, $(b,b,1)$, or $(1)$, and $\varepsilon(\Rot(G))=\varepsilon(G)$ otherwise.
\end{lemma}

\begin{proof} Without loss of generality, we shall  restrict to white-sided graphs. A white-sided graph with the left border vertex being white (resp. black) will be referred to as a white-white (resp. black-white) graph.

The total number of real vertices is odd for white-white graphs and even for black-white graphs. Let $q$ be the number of one-valent real vertices, then 
\begin{equation} \label{eq:polechange} 
  \pol(\Rot(G)) - \pol(G) \equiv 
	\begin{cases}
	  q + 1 \mod 2,&  \text{if~} G \text{ is a white-white graph}, \\
		q \mod 2,&  \text{if~}G \text{ is a black-white graph}. 
	\end{cases}
\end{equation}
This is true since each real vertex of $G$ with $\deg(v) \neq 1$ contributes either to $\pol(G)$ or, after rotation, to $\pol(\Rot(G))$. 

For the case of white-white graphs, the sequences $\Sigma^{\odd}_w$ and $\Sigma^{\odd}_b$ have both odd length, and thus Lemma~\ref{lem_rotationsignchange} implies $\lev(\Rot(G)) \equiv \lev(G) + p \equiv \lev(G) + 1$.
Here $p$ is the total number of even real vertices. 
Together with Equation \eqref{eq:polechange}, we get $\varepsilon(\Rot(G)) = (-1)^q \varepsilon(G)$. 
Hence the sign changes if and only if $\Sigma^{\odd}_w$ is of type $(1,b,b)$, $(b,b,1)$, or $(1)$.

In the case of black-white graphs, each of $\Sigma^{\odd}_w$ and $\Sigma^{\odd}_b$ has even length. Both, only one or none of them can be of type $(1,a)$ or $(a,1)$, corresponding to $q=2,1$ or $0$.
Lemma~\ref{lem_rotationsignchange} implies $\lev(\Rot(G)) \equiv \lev(G) + q$. Together with Equation \eqref{eq:polechange} this gives $\varepsilon(\Rot(G)) = (-1)^{2q}\varepsilon(G)  = \varepsilon(G)$, which proves the claim.
\end{proof}

\begin{remark}
  A simple parity check shows the following: 
	If $d$ is even, all graphs in $\bm{G}$ have border vertices of opposite colour. 
	If $d$ is odd, all graphs in $\bm{G}$ have border vertices of the same colour.  
\end{remark}

\begin{proof}[Proof of Theorem~\ref{Th:InvThGraphs} for $d$ even]
  Rotation  defines a bijection $\Rot : \bm{NW} \to \bm{NB}$ with
	$\varepsilon(G)=\varepsilon(\Rot(G))$ by Lemma~\ref{L:RotSafeGraph}.
	Then Proposition \ref{pro_symmetrizinggraphs} finishes the proof.
\end{proof}

\begin{definition}\label{Def:Reduced}
A \emph{reduced} graph $G$ is a nearly symmetric graph such that $\Sigma_w^{\odd}$ and $\Sigma_b^{\odd}$ are of the form $(a), (1,a,1),$ or $(a,b,b)$ for $a,b \neq 1$.
We denote by $\bm{RW}\subset\bm{N\!W}$ the subset of reduced white-sided graphs.
\end{definition}

\begin{corollary}\label{Cor:Reduced}
  If $d$ is odd, we have $S(\bm{W})=S(\bm{RW})$.
\end{corollary}

\begin{proof}
  Rotation defines an involution on $\bm{NW} \setminus \bm{RW}$ which satisfies 
	$\varepsilon(G)=-\varepsilon(\Rot(G))$ by Lemma~\ref{L:RotSafeGraph}.
	Using Proposition \ref{pro_symmetrizinggraphs}
	we get $S(\bm{RW})=S(\bm{NW})=S(\bm{W})$.
\end{proof}

\subsection{Midline operation and rest of the proof}

In this subsection we assume $d$ odd and set $\bm{RG} = \bm{RW} \cup \bm{RB}$.
We want to construct a bijective correspondence 
$$\Tr:~\{\text{\emph{short} graphs of}~\bm{RG}\}\longrightarrow\{\text{\emph{long} graphs of}~\bm{RG}\}$$ 
which does \emph{not} necessarily respect the colour of a graph.

\begin{figure}[h]
\centering
\includegraphics[scale=0.9]{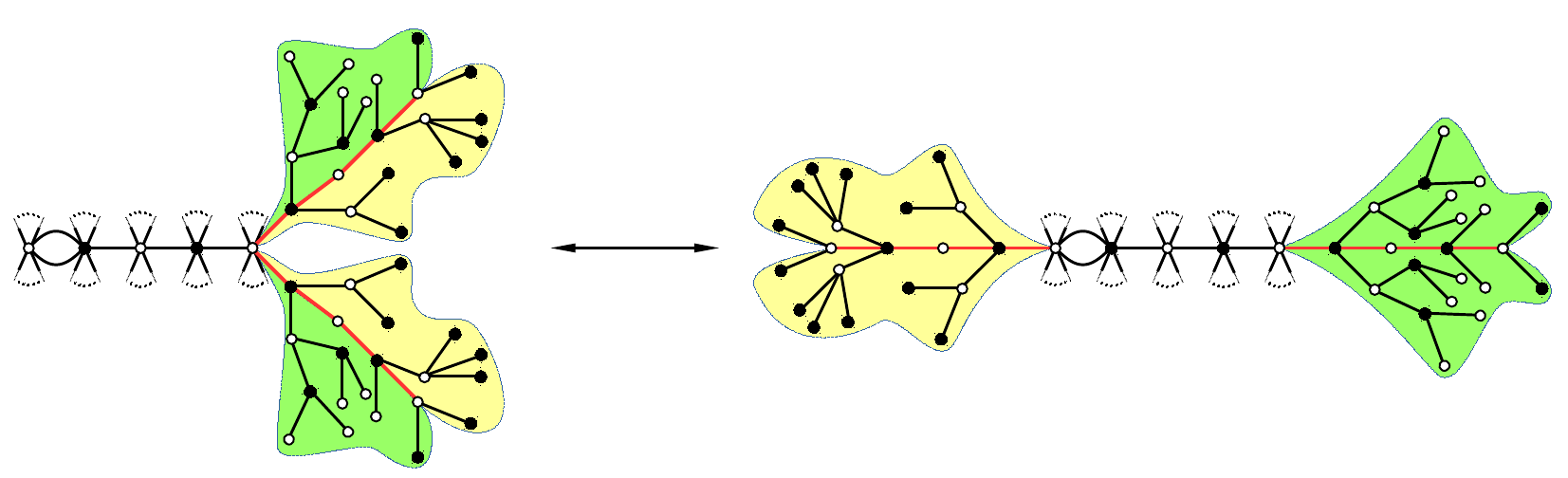}
\caption{The midline cutting method.}
\label{Fig:FirstStep}
\end{figure}

\begin{construction}
Without loss of generality, consider a white-sided short graph $G=G_0\in\bm{RW}$ 
and assume that the cycle lies on the left hand side of $G$ (see bottom graph of Figure~\ref{Fig:TwoGraphExamples}).

\begin{enumerate}
	\item 
	We use the midline cutting method from \cite[Proof of Lemma 3.13]{IZ18}:
	Take the pair of trees growing at the right border vertex closest to the positive
	direction of the real axis and cut each tree into two halves along the
	midline	(see Figure~\ref{Fig:FirstStep}).
	Glue the two $\conj$-symmetric halves closest to the real axis
	to obtain a rooted tree which we attach to the left border vertex of $G$
	such that the glued midline is mapped to $\RR$. Do the same with the second 
	pair of half-trees, but attaching them to the right border vertex. 
	The long graph obtained by this construction is denoted by $G_1$.

	\item 
	The graph $G_1$ is in general not of type $\Lambda_w,\Lambda_b$, since the former border vertices $v_1, v_2$ changed
	their degree. We repair this by applying $G_2:=\Flip_{v_1}^{v_2}(G_1)$. 
	
	\item 
	The graph $G_2$ is in general not nearly symmetric (nor reduced). 
	Indeed, $\Sigma^{\even}_b(G_2)$ is of the form
	$(\alpha,\Xi,\alpha^{-1})$ or $(\alpha,e_0,\Xi,\alpha^{-1})$
	depending on whether $\Sigma^{\even}_b(G) = \Xi$ or 
	$\Sigma^{\even}_b(G) = (e_0,\Xi)$ (here and in the following, 
	$\Xi$ denotes a symmetric sequence). 
	Similarly, $\Sigma^{\even}_w(G_2)$ is of the form
	$(\alpha,\Xi,e_0,\alpha^{-1})$ or $(\alpha,\Xi,e_1,e_1,\alpha^{-1})$
	depending on whether $\Sigma^{\even}_w(G) = (e_0,\Xi)$ or 
	$\Sigma^{\even}_w(G) = (e_1,\Xi,e_1)$. 
	The first case is symmetric, while in the latter three cases 
	we perform a cyclic shift on the subsequences 
	$(\Xi,e_1,e_1)$, $(\alpha,e_0)$, and $(\alpha,\Xi,e_0)$,
	respectively, to make the sequences nearly symmetric. 
	Finally, if either $\Sigma^{\odd}_b(G_2)$ or $\Sigma^{\odd}_w(G_2)$
	is of type $(b,a,b), a,b\neq 1$, we apply a flip on $(b,a)$.
	We obtain a reduced nearly symmetric long graph $G_3$, 
	which we also denote $\Tr(G)$.  
 \end{enumerate}
\end{construction}
 
\begin{lemma}\label{L:BijTr}
  The map $\Tr$ is a bijection.
\end{lemma} 

\begin{proof}
  Given a reduced nearly symmetric extended graph, the position of the cycle indicates the 
	length of the two real segments on the left and right side of the graph that we want to reorganize. 
	In particular, the length of $\alpha$ in step (3) can be recovered from $G$ and hence, knowing 
	the real sequences of $G$, step (3) can be reversed uniquely. 
	Next, we can undo step (2) by flipping the two vertices which are at the 
	inner ends of these two segments. 
	Finally, there is unique way of cutting open the two segments and 
	regluing them as additional tree at the right border vertex to undo step (1).
\end{proof}

\begin{lemma}\label{L:SignTr}
Let $G \in \bm{RW}$ be a short graph. Then 
 
\begin{equation} \label{eq:TransformSign} 
  \varepsilon(\Tr(G)) = 
	\begin{cases}
	  - \varepsilon(G) & \text{ if } \Tr(G) \text{ is white-sided}, \\
		+ \varepsilon(G) & \text{ if } \Tr(G) \text{ is black-sided}. 
	\end{cases}
\end{equation}
\end{lemma}

\begin{proof}
  Throughout the proof, $\Sigma$ and $\overline{\Sigma}$ will refer to sequences for $G$ and $\Tr(G)$ respectively. By symmetry we may restrict our attention to graphs $G$ with the cycle sitting on the left hand side. Fix a colour $c$ and let $k_c$ be the number of transpositions needed to transform $\overline{\Sigma}_c$ to $(\overline{\Sigma}_c^{\odd}, \overline{\Sigma}_c^{\even})$, as in the proof of Lemma \ref{lem_rotationsignchange}. Set $l_c = l(\overline{\Sigma}_c^{\even})$ and let $n_c$ be the number of even vertices of colour $c$ located \emph{before} the cycle in $\Tr(G)$. In other words, $n_c = l(\alpha)$ for the sequence $\alpha$ appearing in step (3). Note that the cyclic shifts applied in step (3) do not affect the number $k_c$ and we can express it as
\begin{equation} \label{eq:disDiff1} 
  k_c = 
	\begin{cases}
	  n_c & \text{ if } l(\overline{\Sigma}_c^{\odd}) = 1, \\
		n_c + l_c & \text{ if } l(\overline{\Sigma}_c^{\odd}) = 3. 
	\end{cases}
\end{equation}
	Let $\delta_c$ denote the \emph{change} of pole disorders 
	from $G$ to $\Tr(G)$ caused by vertices of colour $c$. We have $\delta_c = n_c$ for $\overline{\Sigma}_c^{\odd} = (a)$ or $(1,a,1)$ and $\delta_c = n_c + 1$ for $\overline{\Sigma}_c^{\odd} = (a,b,b)$.
	Hence the total contribution of the colour $c$ to the sign change can be expressed as
\begin{equation} \label{eq:disDiff2} 
  \dis(\overline{\Sigma}_c) - \dis((\Sigma_c^{\odd}, \Sigma_c^{\even})) + \delta_c  \equiv 
	\begin{cases}
	  0 & \text{ if } l(\overline{\Sigma}_c^{\odd}) = 1, \\
		l_c + 1 & \text{ if } l(\overline{\Sigma}_c^{\odd}) = 3. 
	\end{cases}
\end{equation}

Since by assumption $G$ is a white-sided short graph with the cycle located to the left-hand side, we have $\Sigma_b = (\Sigma_b^{\odd}, \Sigma_b^{\even})$ and $\Sigma_w = (\Sigma_w^{\even}, \Sigma_w^{\odd})$. Hence for case $c=w$ we should take into account an extra contribution of $l_w$, and thus the sign change for black and white vertices is given in the following table.

	\begin{center}
	\begin{tabular}{ccc}
\hline
  & $w$ & $b$ \\ \hline
$l(\overline{\Sigma}_c^{\odd}) = 1$ & $(-1)^{l_w}$   & $+1$ \\
$l(\overline{\Sigma}_c^{\odd}) = 3$ & $-1$   & $(-1)^{l_b + 1}$ \\ \hline
\end{tabular}
	\end{center}
	
	In the case where $\Tr(G)$ is black-sided, we multiply the diagonal entries of the table.
	Since $l_w + l_b$ is odd, we obtain the total sign change $+1$. In the case
	where $\Tr(G)$ is white-sided, we multiply the antidiagonal entries and obtain $-1$. 
\end{proof}

\begin{proof}[Proof of Theorem~\ref{Th:InvThGraphs} for $d$ odd]
Following Corollary~\ref{Cor:Reduced}, it is sufficient to prove that $S(\bm{RW})=S(\bm{RB})$. 
The latter equality follows from Lemmata~\ref{L:BijTr} and~\ref{L:SignTr}.
\end{proof}

%% file: invariance.tex
\section{Invariance theorem} \label{sec_invariance}

In this section we prove Theorem~\ref{Th:InvarTh}. 
We fix $d$ and partitions $\Fatlambda = (\Lambda_1,\ldots,\Lambda_k)$ of $d$ such that
$d = \sum_{i,j} (\Lambda_j^i - 1)$.
First note that the number of functions contained in $\SSS(\PPP, \Fatlambda)$ only depends on the ordering
of the branch points in $\PPP=(b_1, \dots, b_k)$ on the real line, not on their exact position. 
This follows from Theorem \ref{thm_rationalfunctionasdessins}.
Hence it remains to prove that the signed count does not change under a permutation of the points in $\PPP$. 
Based on Section \ref{sec_dessins} it is slightly more convenient to take the symmetric point of view here:
We fix the ordering of points by $b_1 < \dots < b_k$, and instead permute the partitions
stored in the sequence $\Fatlambda$. Of course, it suffices to prove invariance
under the action of the transposition exchanging  $\Lambda_j$ and $\Lambda_{j+1}$, $j=1, \dots, k-1$.

Let us fix some $j \in \{1, \dots, k-1\}$. The relation to the $bw$-graphs from the previous section is as follows.
We rename $\Lambda_j=\Lambda_b$ and $\Lambda_{j+1}=\Lambda_w$. Let $\Gamma$ be a real simple rational dessin of degree $d$ and type $\Fatlambda$. We colour the preimages of $b_j$ and $b_{j+1}$ in black and white, respectively. 
Let $\BWG_{bw} := \varphi^{-1}([b_j, b_{j+1}]) \subset\Gamma$ denote the preimage of the closed segment joining $b_j$ and $b_{j+1}$. 
We need to distinguish two cases.
Let us recall from \cite[Definition 3.1]{IZ18} that a \emph{black and white tree}, or $bw$-tree, is a graph as described in Definitions~\ref{Def:BWGraphs},
except for requiring a tree graph, in exchange for our cycle condition.

\begin{figure}
\centering
\includegraphics[scale=0.9]{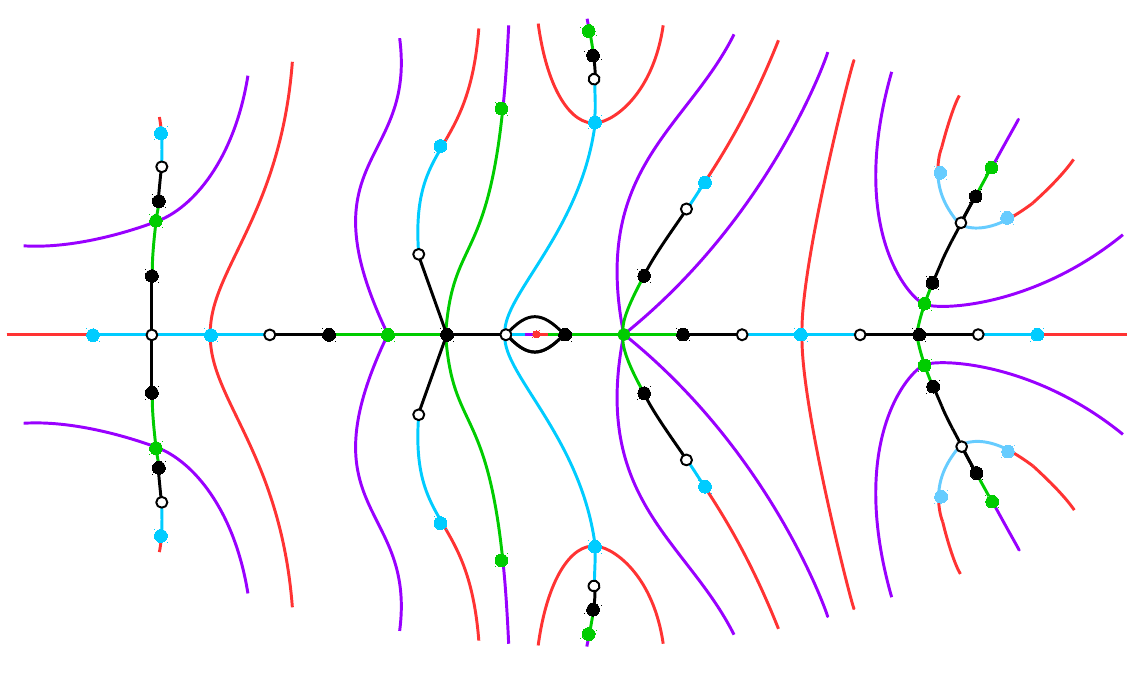}
\caption{An example of a real dessin satisfying condition \ref{bw2}.}
\label{Fig:GraphDessin}
\end{figure}

\begin{enumerate}[label=(\arabic*)]
	\item \label{bw1} $\BWG_{bw}$ is a union $T_1 \cup \ldots \cup T_r$ of disjoint $bw$-trees, 
	where some of them are real and others split into pairs of trees that are complex conjugate.

	\item \label{bw2} $\BWG_{bw}$ is a union $R_0 \cup T_1 \cup \ldots \cup T_r$ of $bw$-trees $T_1,\ldots,T_s$ 
				as before, together with exactly one real simple $bw$-graph $R_0$
	      (see Figure~\ref{Fig:GraphDessin}).
\end{enumerate}

We focus on case \ref{bw2}. Case \ref{bw1} is similar to the situation in \cite{IZ18} and will be treated later.
The following definition describes the type of graphs we obtain when we make $b_j$ and $b_{j+1}$ \enquote{collide} in $\RR$,
or equivalently, collapse the subgraph $\BWG_{bw} \subset \Gamma$. 
A \emph{real polynomial dessin} $\tilde{\Gamma} \subset \CC\PP^1$ is a labelled graph satisfying all
the conditions from Definition \ref{Def:RealSimpleDessin} except for condition (c), which we replace by
\begin{enumerate}
	\item[(c')] exactly one vertex is labelled by $\infty$, namely $\infty \in \CC\PP^1$. 
\end{enumerate}

\begin{definition}\label{Def:iEnhDess}
A \emph{$bw$-enhanced dessin} is a real polynomial dessin $\tilde{\Gamma}$ with (ordered) labels $\{1, \dots,j-1, bw, j+2, \dots, \infty\}$
together with the following data:
\begin{itemize}
  \item A real $bw$-vertex $V_0$, together with two partitions $\Pi_b(V_0),\Pi_w(V_0)$ of $\deg(V_0)/2+1$.
	\item Let $\mathcal{V}_{\R}$ be the set of real $bw$-vertices different from $V_0$ and let $\mathcal{V}_+$ be the set of $bw$-vertices with positive imaginary part.
	      For each vertex $V \in \mathcal{V}_{\R} \cup \mathcal{V}_+$, fix two partitions $\Pi_b(V),\Pi_w(V)$ of $\deg(V)/2$.

\end{itemize} 

\end{definition}

From a simple rational dessin $\Gamma$ of type \ref{bw2} 
we can construct a $bw$-enhanced real dessin $\tilde{\Gamma}$ as follows. 
First, we obtain the graph $\tilde{\Gamma}$ by contracting all edges $j\rightarrow j+1$ of $\Gamma$. 
This means that each component of $\BWG_{bw} = R_0 \cup T_1 \cup \ldots \cup T_r$
contracts to a point, and we label these points by $bw$. 
We set $V_0$ the vertex of $\tilde{\Gamma}$ obtained from contracting $R_0$ and define
$\Pi_b(V_0)$ and $\Pi_w(V_0)$ as the partition of degrees of black and white vertices in $R_0$, respectively.
Similarly, for each vertex $V \in \mathcal{V}_{\R} \cup \mathcal{V}_+$, 
we define $\Pi_b(V_0)$ and $\Pi_w(V_0)$ as the partition of degrees of black and white vertices in $T_i$, where $T_i$
is the $bw$-tree contracting to $V$.
We call $\tilde{\Gamma}$ the \emph{contraction} of $\Gamma$.

\begin{definition}\label{Def:SpecDis}
Let $\tilde{\Gamma}$ be a $bw$-enhanced dessin. For any label $j \neq bw, \infty$, we set
$\Sigma_j$ to the sequence of degrees of real vertices labelled by $j$, ordered from left to right.
We denote by $\tilde{\Sigma}_w$ and $\tilde{\Sigma}_b$ the sequences of partitions $\Pi_w(V)$ and $\Pi_b(V)$,
respectively, running through the real $bw$-vertices (from left to right), except for $V_0$, for which 
we use the partitions $(\Pi_w(V),1)$ and $(\Pi_b(V),1)$. A disorder of $\tilde{\Sigma}_w$ or $\tilde{\Sigma}_b$
is a pair of entries $\Pi_w^k(V) > \Pi_w^{k'}(V')$ or $\Pi_b^k(V) > \Pi_b^{k'}(V')$, respectively, with $V < V'$. 
Clearly, multiple entries of a partition are counted multiple times here. 
We set
\[
  \dis(\tilde{\Gamma}) := \dis(\tilde{\Sigma}_w) + \dis(\tilde{\Sigma}_b) + \sum_{j \neq bw, \infty} \dis(\Sigma_j), 
	\quad \quad \quad \quad
	\varepsilon(\tilde{\Gamma}) := (-1)^{\dis(\tilde{\Gamma})}.
\]

\end{definition}

Recall from~\cite{IZ18} that the sign $\varepsilon(T)$ of a real $bw$-tree is $(-1)^{\dis(\Sigma_w) + \dis(\Sigma_b)}$.

\begin{lemma}\label{L:SignEnh=SignGamm}
If a $bw$-enhanced dessin $\tilde{\Gamma}$ is the contraction of a dessin $\Gamma$, then 
\[
  \varepsilon(\Gamma)=\varepsilon (\tilde{\Gamma})\cdot\varepsilon(R_0)\cdot\prod_T\varepsilon (T),
\]
where the product is taken over the real trees in $\BWG_{bw}$.
\end{lemma}

\begin{proof}
  The right hand side of the formula takes care of all disorders contributing to $\varepsilon(\Gamma)$ except 
	for the ones involving a vertex in the interior of the (topological) edge $E$ of $\Gamma$ containing $p$. 
	Note that for each label $i \in \{j+2, \dots k, \infty, 1, \dots, j-1\}$ there is exactly one vertex
	labelled by $i$ in the interior of $E$, and all these vertices have valance to $2$. Hence such disorders
	appear in pairs: Any higher-degree real vertex to the left of $p$ produces exactly one level and one pole disorder.
	This proves the formula. 
\end{proof}

Let $\tilde{\Gamma}$ be a $bw$-enhanced dessin. 
For any real $bw$-vertex $V$, we define the colour $c(V)$ to be white if the real edge to the right of $V$ is of type $bw \to j+2$. Otherwise we set
$c(V)$ to be black.
We denote by $\mathcal{R}_0$ the set of $c(V)$-sided real simple $bw$-graphs of type $\Pi_b(V_0),\Pi_w(V_0)$.
For any $V \in \mathcal{V}_\R$ we denote by $\mathcal{T}_V$ the set of $c(V)$-sided real $bw$-trees of type $\Pi_b(V),\Pi_w(V)$.
For any $V \in \mathcal{V}_+$ we denote by $\mathcal{T}_V$ the set of (non-real) $bw$-trees of type 
$\Pi_b(V),\Pi_w(V)$ with one marked half-edge emanating from a white vertex.

For any $V \in \mathcal{V}_+$, let us mark an adjacent edge $_V$ of type $bw \to j+2$.
Let $\mathcal{M}$ denote the set of increasing real simple dessins $\Gamma$ that contract to $\tilde{\Gamma}$. Then we can define a map
\begin{align*} 
  F : \mathcal{M} &\to \mathcal{R}_0\times\prod_{V\in\mathcal{V}_\R\cup\mathcal{V}_+}\mathcal{T}_V, \\
	    \Gamma     &\mapsto (R_0, T_V : V\in\mathcal{V}_\R\cup\mathcal{V}_+),
\end{align*}
where $R_0$ and $T_V$ denote the components of $\BWG_{bw}$ contracting to $V_0$ and $V$, respectively. 
For $V \in \mathcal{V}_+$, we additionally mark the first half-edge in $T_V$ 
touching $e_V$, counted in clockwise direction.

\begin{lemma}\label{L:HatGammaAndGamma}
  The map $F$ is a bijection.
\end{lemma}

\begin{figure}
\centering
\includegraphics[scale=0.8]{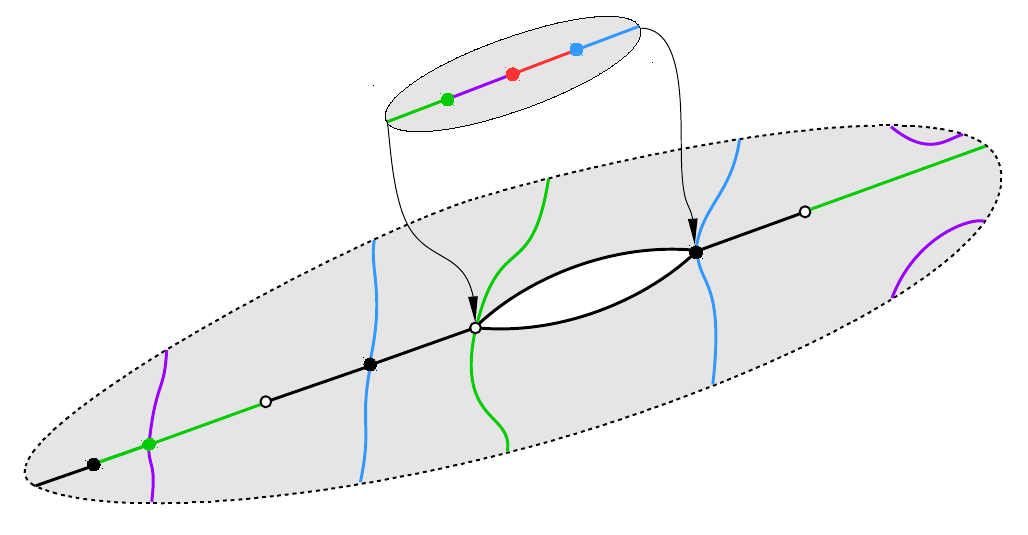}
\caption{Filling the \enquote{hole} of $R_0$.}
\label{Fig:GluingD1}
\end{figure}

\begin{proof}
  We can easily describe the inverse map. Given $(R_0,~T_V : V\in\mathcal{V}_\R\cup\mathcal{V}_+)$, we insert $R_0$ and $T_V$ into a small disc removed 
  around $V_0$ and $V \in \tilde{\Gamma}$, respectively, as described in~\cite[Proof of Lemma 4.6]{IZ18}.
	Additionally, we fill the ``hole'' bounded by the cycle in $R_0$ as follows. 
	Let $D \subset \CC$ denote the unique affine real polynomial dessin of degree $1$ labelled by 
	$-\infty \to j+2 \to \dots \to k \to \infty \to 1 \to \dots \to j-1 \to +\infty$. 
	We glue $D$ into the hole such that the edge $-\infty \to j+2$ is attached to the white 
	vertex in the boundary (see\ Figure~\ref{Fig:GluingD1}). 
	In this way we obtain real simple rational dessin $\Gamma$ such that $F(\Gamma)$ is equal to the given data. 
\end{proof}

\begin{corollary} \label{L:SnumProd}
  With the notation used in Lemma \ref{L:HatGammaAndGamma}, we have
	\[
	  S(\mathcal{M}) = \varepsilon(\tilde{\Gamma}) \; S(\mathcal{R}_0) \prod_{ V\in\mathcal{V}_\R} S(\mathcal{T}_V) \prod_{ V\in\mathcal{V}_+} |\mathcal{T}_V|.
	\]
\end{corollary}

\begin{proof}
  Follows from Lemma \ref{L:SignEnh=SignGamm} and Lemma \ref{L:HatGammaAndGamma}.
\end{proof}

\begin{proof}[Proof of Theorem~\ref{Th:InvarTh}]
  Since transpositions of the form $(j, j+1)$ generate the symmetric group, it is enough
	to prove $S(\SSS(\PPP, \Fatlambda)) = S(\SSS(\PPP, \Fatlambda'))$, 
	where $\Fatlambda'$ is the sequence of partitions with $\Lambda_j$ and $\Lambda_{j+1}$ swapped.
	Let $\SSS^{\text{\ref{bw1}}}(\PPP, \Fatlambda)$ and $\SSS^{\text{\ref{bw2}}}(\PPP, \Fatlambda)$ 
	be the subsets of functions/dessins of type \ref{bw1} and \ref{bw2}, respectively.
	Then $S(\SSS^{\text{\ref{bw1}}}(\PPP, \Fatlambda)) = S(\SSS^{\text{\ref{bw1}}}(\PPP, \Fatlambda'))$
	follows, after some straightforward adaptation to our case, from \cite[Lemma 4.8]{IZ18}.
  It remains to show $S(\SSS^{\text{\ref{bw2}}}(\PPP, \Fatlambda)) = S(\SSS^{\text{\ref{bw2}}}(\PPP, \Fatlambda'))$.
	We can prove this for each contraction separately. More precisely,
	let $\tilde{\Gamma}, \tilde{\Gamma}'$ be a pair of $bw$-enhanced dessins obtained from each other by flipping the 
	partitions $\Pi_b(V), \Pi_w(V)$ for any $bw$-vertex $V$. Let $\mathcal{M}$ and $\mathcal{M}'$ 
	be the sets of dessins of type $\Fatlambda$ and $\Fatlambda'$ that contract to 
	$\tilde{\Gamma}$ and $\tilde{\Gamma}'$, respectively. We show 
	$S(\mathcal{M})=S(\mathcal{M}')$ using Corollary \ref{L:SnumProd}.
	Indeed, the invariance of the various factors follows for $\varepsilon(\tilde{\Gamma})$ by definition, 
	for $S(\mathcal{R}_0)$ by Theorem \ref{Th:InvThGraphs}, for $S(\mathcal{T}_V)$ by 
	\cite[Theorem 9]{IZ18}, and for $|\mathcal{T}_V|$ by switching the colours of the vertices.
	Note that in the last case, the number of half-edges adjacent to a white and black vertex, respectively, is equal
	and hence the number of possible markings does not change.
\end{proof}

\begin{figure}
		\centering
		\includegraphics[scale=1]{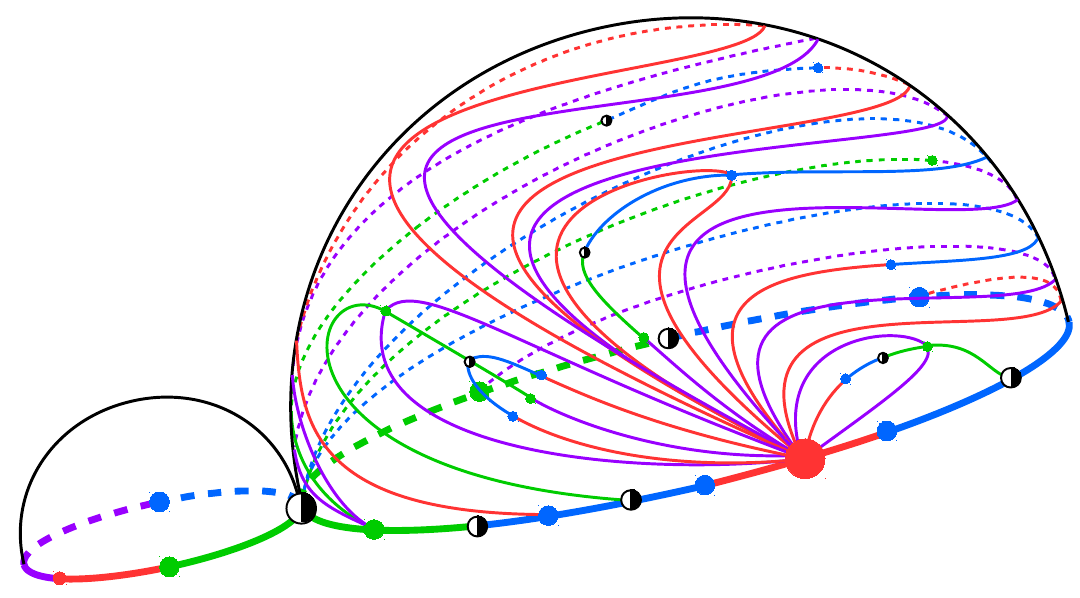}
		\caption{The bubbling of a dessin into a reducible dessin of degree $(d-1, 1)$.}
		\label{Fig:Bubbling}
	\end{figure}

\begin{remark}
  Despite the combinatorial nature of the presented proof of Theorem~\ref{Th:InvarTh}
	it is instructive to also get a geometric picture of the degeneration when $b_j$ and $b_{j+1}$ collide.
  While in our description the two discs adjacent to $p$ just collapse to $V_0$ in $\tilde{\Gamma}$,
	the limit of this degeneration in $\overline{\MMM}_{0,0}(\CC\PP^1, d)$ is a ramified cover with 
	reducible source curve. The \enquote{bubbling} produces a source curve with two components
	(see Figure~\ref{Fig:Bubbling}).
	The main component is described by $\tilde{\Gamma}$ and is mapped to $\CC\PP^1$ by a polynomial map 
	of degree $d-1$. The bubble corresponds to the two discs adjacent to $p$ and is mapped 
	isomorphically to $\CC\PP^1$. The vertex $V_0$ is the node connecting the two spheres.
	This phenomenon does not occur in the polynomial case of \cite{IZ18}.
	What keeps the situation under control for simple rational functions is the fact
	that the degree $1$ map on the bubble is essentially unique and hence can be discarded
	from the combinatorial data.
\end{remark}

%% file: vanishing.tex
\section{Vanishing Statements} \label{sec vanishing}

This section is devoted to proving the vanishing statements contained in Theorem \ref{Th:IffOdd} or, in other words,
the \emph{only if} part.
Let $\bm{\lambda}:=(\lambda_1,\ldots,\lambda_k)$ be a sequence of \emph{reduced} partitions
as in the introduction.

\begin{proposition}[First part of non-vanishing for all degrees]\label{Th:OnlyIfEven}
  If a partition $\lambda_j$ has more than one odd number appearing an odd number of times or more than one even number appearing an odd number of times, 
	then $F_{\bm{\lambda}}^\odd(q) = 0$ and $F_{\bm{\lambda}}^\even(q) = 0$.
\end{proposition}

\begin{proof}
  Fix $a\neq b$ of the same parity which both appear an odd number of times in $\lambda_j$.
	Then for any dessin $\Gamma$ there is an odd number of real vertices 
	of label $j$ and degree $2a+2$ and $2b+2$, respectively.
	Let $v_1 < \dots < v_{2m}$ be the totality of such vertices. 
	We want to invert the order of appearance of these vertices by applying a flip operation
	of the form
	\begin{equation} \label{eq:flipping} 
	  \Flip^{v_1}_{v_{2m}} \circ \dots \circ \Flip^{v_m}_{v_{m+1}}.
	\end{equation}
	where $\Flip^{v}_{w}$ is the following adaptation of Definition \ref{Def:flip}.
	\begin{enumerate}
		\item 
		  If $v$ and $w$ have the same degree, we do nothing. 
		\item 
		  If $v$ and $w$ have degree $2a+2$ and $2b+2$, respectively, and $a > b$,
			we set $F(v) := (T_1, \dots, T_{2(a-b)})$ to be the sequence of pairs of 
			trees growing at $v$ and closest to the negative direction of $\RR$. 
			If $a < b$, we proceed symmetrically.
		\item 
		  The real segments near $v,w$ can be of the following types.
			\begin{align*} 
				a,b \text{ even:} & & j-1\rightarrow j \rightarrow j+1 & & \text{ or } & & j+1\rightarrow j \rightarrow j-1 \\
				a,b \text{ odd:} & & j-1\rightarrow j \rightarrow j-1 & & \text{ or } & & j+1\rightarrow j \rightarrow j+1 
			\end{align*}
			If $v$ and $w$ are of the same type, we remove $F(v)$ from $v$  and glue it to $w$
			(say, closest to the negative direction of $\RR$, but not in the interior of the cycle, if present). If the vertices are of opposite type, we glue back the inverted sequence $F(v)^{-1}$
			instead. In this way, we ensure that the pattern of increasing and decreasing edges alternates when going around $w$, as required. 
	\end{enumerate}
	In summary, we see that operation \eqref{eq:flipping} defines an involution on the
	set of dessins which inverts the order of appearance of the vertices 
	$v_1, \dots, v_{2m}$.
	It is shown in \cite[Page 37]{IZ18} that this operation 
	flips the sign of a dessin. Note that cited argument can be applied 
	to our situation since the number of pole disorders is not affected by \eqref{eq:flipping}. 
	Hence the statement follows.
%
	%
	%
	%
\end{proof}
	
\begin{proposition}[non-vanishing for odd degrees]\label{Prop:OnlyIfOdd}
Consider a sequence of partitions $\bm{\lambda}$ which does \emph{not} satisfy the conditions of Proposition~\ref{Th:OnlyIfEven}. 
If there is an odd number of partitions $\lambda_j$ having an even element appearing an odd number of times, then $$F^{\odd}_{\bm{\lambda}}(q)= 0.$$
\end{proposition}

\begin{proof}
For $d$ odd the operation $\varphi(z) \mapsto \tilde{\varphi}(z) = \varphi(-z)$ produces an increasing function again and hence
defines an involution on $\SSS(\PPP,\Fatlambda)$. 
Let $\Sigma=\Sigma_{b_j}$ the ramification sequence of a critical value as in Definition \ref{def_signsimplerationalfunction}.
Then the corresponding sequence for $\tilde{\varphi}$ is given by $\Sigma^{-1}$, the inverted sequence.
One easily checks
$$\dis(\Sigma) + \dis(\Sigma^{-1})\equiv_{2} c(\Sigma),$$ 
where 
$$c(\Sigma):=\sharp\{\text{pairs of numbers}~ m, n~\text{in}~\Sigma~\text{such that}~m\neq n\}.$$ 
Recall that $\Sigma$ includes a $1$ coming from the simple pole, so the sum of its entries is even.
Now, if $\lambda_i$ contains an even number $a$ an odd number of times, 
then $\Sigma$ contains two or three elements appearing an odd number of times ($1$, $a+1$, and possibly an even entry), 
and hence $c(\Sigma)$ is odd. 
If, on the other hand, $\lambda_j$ does not contain an even element an odd number of times, 
then $\Sigma$ contains zero or one elements appearing an odd number of times
and hence $c(\Sigma)$ is even. 
The result then follows by the previous formula.
\end{proof}

\begin{remark} 
  For $d$ even, we could instead try to use the transformation $\varphi(z) \mapsto -\varphi(-z)$,
	which provides a bijection between $\SSS(\PPP,\Fatlambda)$ and 
	$\SSS(-\PPP,\Fatlambda)$. The same argument as above shows that 
	$F^{\even}_\fatlambda(q) = 0$ if there is an odd number of partitions $\lambda_j$
	with an odd entry appearing an odd number of times. 
	But note that such entries correspond 
	to local extrema of $\varphi$. Hence the number is even.
\end{remark}

%% file: brokensequences.tex
\section{Broken alternations}\label{Sec:BrokenAlt}

\begin{definition}
  Let $a \in \Sym(n)$ be a permutation on $n$ elements, $i \mapsto a_i$. 
	We call $a$ an \emph{(ordinary) alternation} if 
	$a_n < a_{n-1} > a_{n-2} < a_{n-3} > \dots \lessgtr a_1$.  
	The permutation $a$ is said to be a \emph{broken alternation} if it is an alternation
	except for exactly one index $i \in \{1, \dots, n-1\}$, where
	$a_i,a_{i+1}$ is increasing if it should be decreasing,
	or the opposite. In this case, $i$ is called the \emph{break} of $a$, and we write $(a_1 \dots a_i | a_{i+1} \dots a_n)$. 
	We denote by $A_n$, $B_n$ and $B_n^j$ the number 
	of ordinary alternations, broken alternations, broken alternations with $a_j = n$, respectively.
\end{definition}

It is clear that we have $A_1 = A_2 = B_2 =1$, $B_1=0$, and $B_n= B_n^1 + \dots + B_n^n$.

\begin{proposition}\label{Prop:AlterDessins}
  There is a one-to-one correspondence between broken alternations of $n$ elements 
	and real simple rational dessins of degree $n$ and type $\Fatlambda = (\Lambda_1, \dots, \Lambda_n)$, 
	$\Lambda_1 = \dots = \Lambda _n = (2 1 \dots 1)$. 
	Moreover, the sign of any such real simple rational dessin is equal to $-1$ to the power of $\left\lfloor n/2 \right\rfloor$.
\end{proposition}

\begin{figure}
\centering
\includegraphics[scale=2]{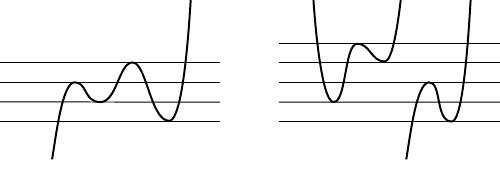}
\put(-162,25){\small{$\cdots 1\cdots$}}
\put(-162,36.6){\small{$\cdots 2\cdots$}}
\put(-162,48.2){\small{$\cdots 3\cdots$}}
\put(-162,59.8){\small{$\cdots 4\cdots$}}
\put(-148,71.4){\small{$ 5\cdots$}}
\caption{To the left, a polynomial representing the ordinary alternation $a=(3241)$. To the right, a simple rational function representing the broken alternation $a=(254|31)$.}
\label{Fig:ExampleAlternations}
\end{figure}
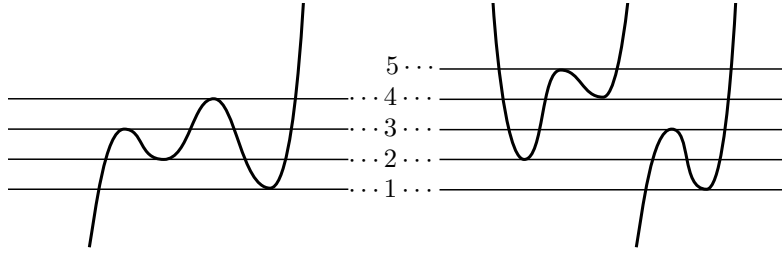

\begin{remark} \label{rem_ordinaryalternations}
  The exact same relation holds between ordinary alternations and real polynomial dessins,
	see \cite[Example 1.5]{IZ18}. We refer to Figure \ref{Fig:ExampleAlternations} for examples.
\end{remark}

\begin{proof}[Proof of Proposition~\ref{Prop:AlterDessins}]
	Given a real simple rational dessins of type $\Fatlambda$, the
	sequence of labels of the real vertices of degree greater than $2$,
	ordered from left to right, 
	produces a broken alternation of length $n$.
	The break corresponds to passing through the simple pole $p$.
	Vice versa, given a broken alternation $a$, we can construct a dessin as follows.
	We place $4$-valent vertices on $\R$ with labels as described in $a$. 
	We put a $2$-valent vertex with label $\infty$ on the segment between the vertices
	corresponding to the break of $a$. We glue the imaginary ends of the graph to
	infinity, except for the ends to the left and right of $p$, which we glue
	pairwise to obtain a pair of cycles. Finally we insert $2$-valent vertices 
	in the unique way to complete the graph to a dessin.
	The sign is deduced by noting that an real simple rational function of type $\Lambda _n$
	has exactly $\lfloor n/2\rfloor$ local maxima, and each local maxima contributes by a factor of $-1$.
\end{proof}

We can easily derive recursive formulas for the numbers $B_n, B_n^j$ similar to those for $A_n$.
In accordance with Proposition \ref{Prop:AlterDessins}, we set $B_0 = 0$, $A_0 = 1$.

\begin{proposition}\label{Prop:Recursion}
For $n \geq 2$ we have 
\[
 B_n^j =
\begin{cases}
  \binom{n-1}{j-1} \cdot (B_{j-1} A_{n-j} + A_{j-1} B_{n-j}) & \text{if $n+j$ odd,} \\
	A_{n-1} 																									 & \text{if $j=n$ or $n$ odd, $j=1$,} \\
  0 																												 & \text{otherwise.} 
\end{cases}
\]

\end{proposition}

\begin{proof}
By our conventions, the pattern for ordinary alternations is 
$a_{j-1} < a_j > a_{j+1}$ for $n+j$ odd, and $a_{j-1} > a_j < a_{j+1}$ otherwise.
In a broken alteration, we can flip at most one of these inequalities.
Hence, if $a$ is a broken alternation with $a_j = n$, we have the following cases.
If $n+j$ is even, it is of the form $a = (n|a')$ 
or $a = (\overline{a'}|n)$, 
where $a'$ is an ordinary alternation of length $n-1$ 
and $\overline{a'}$ is
defined by $\overline{a'_i} = n - a'_i$.
The map $a \mapsto a'$ is obviously
a bijection, which settles the second and third case. If $n+j$ is odd,
the break of $a$ is different from $j-1$ and $j$. Hence $a$ is of the form
$(a' n a'')$, where either $a'$ is ordinary and $a''$ is broken, or vice versa.
To be precise, this is true after using the unique ordered relabelling of 
the values of $a', a''$ to $\{1, \dots, j-1\}$ and $\{1, \dots, n-j\}$, respectively.
On other hand, a (relabelled) pair $(a',a'')$ occurs exactly 
$\binom{n-1}{j-1}$ times in these constructions, 
since this is the number of ordered maps $\{1, \dots, j-1\} \to \{1, \dots, n-1\}$
(or, equivalently, $\{1, \dots, n-j\} \to \{1, \dots, n-1\}$).
\end{proof}

The first few values of $B_n$ are given in the following table.
\[
	\begin{array}{c|ccccccccccccc}
		n				& 1 & 2 & 3 & 4 & 5  & 6   & 7   & 8    & 9     & 10     & 11      & 12 \\
		\hline
		B_n 		& 0 & 1 & 2 & 7 & 26 & 117 & 594 & 3407 & 21682 & 151853 & 1160026 & 9600567 
	\end{array}
\]	
It seems that the sequence has not appeared in the literature before (cf.\ \cite{Sta})
and is not recorded on \href{https://oeis.org/}{OEIS}.

We can turn the recursive relations into differential equations 
for generating series. 
Since it is more convenient to do this for odd and even indices separately, 
we define for the ordinary alternations
\[
  f(q):=\sum_{k=0}^\infty (-1)^k \frac{A_{2k+1}}{(2k+1)!}q^{2k+1},\quad 
	g(q):=\sum_{k=0}^\infty (-1)^k \frac{A_{2k}}{(2k)!}q^{2k},
\] 
and for the broken alternations 
\[
  u(q):=\sum_{k=0}^\infty (-1)^k \frac{B_{2k+1}}{(2k+1)!}q^{2k+1},\quad 
	v(q):=\sum_{k=0}^\infty (-1)^k \frac{B_{2k}}{(2k)!}q^{2k}. 
\]

\begin{lemma}\label{L:ODE}
  The power series $f,g,u,v$ satisfy the following differential equations.
	\begin{align}
		u' &= 2(g- fu -1), \label{eq:ODE1}\\
		v' &= -(f + fv + gu). \label{eq:ODE2}
	\end{align}
\end{lemma}

\begin{proof} 
	From Proposition~\ref{Prop:Recursion}, for $k \geq 1$ we have 
	\begin{align}\label{eq:B2k+1}
		B_{2k+1} &= 2A_{2k} + \sum_{j=1}^k \binom{2k}{2j-1}\cdot \left( B_{2j-1}A_{2k+1-2j} + A_{2j-1}B_{2k+1-2j} \right) \\
						 &= 2 \left( A_{2k} + \sum_{j=1}^k \binom{2k}{2j-1}\cdot B_{2j-1}A_{2k-(2j-1)} \right).
	\end{align} 
	which proves the first equation for all non-constant terms. For the constant term, we note that by our conventions
	$B_1 = 2(A_0 -0 -1) = 0$. 
	Again by Proposition~\ref{Prop:Recursion}, for $k \geq 1$ we have 
	\begin{align}
	  B_{2k} &= A_{2k-1} + \sum_{j=0}^{k-1} \binom{2k-1}{2j} \cdot \left( B_{2j}A_{2k-2j-1} + A_{2j}B_{2k-2j-1} \right), \\
		       &= A_{2k-1} + \sum_{j=0}^{2k-1} \binom{2k-1}{j} \cdot B_{j}A_{2k-1-j},\label{eq:B2k}
	\end{align}
  The sum in \eqref{eq:B2k} corresponds to the anti-symmetric part of $(f+g)(u+v)$ 
	which is $fv + gu$. 
	Moreover, because of the antisymmetry of the second equation, no constant term appears,
	which completes the proof. 
\end{proof}

We can solve these differential equations explicitly. 
Recall the 19th century result \cite{An81} (see also \cite{Sta}) that 
\[
  f(q) = \tanh(q) \quad \text{and} \quad g(q)=\sech(q). 
\] 
\begin{corollary} \label{uvpolynomialinfg}
  The generating series $u$ and $v$ can be expressed as
	\begin{align}\label{eq:ODE12Sol}
		u(q) &= -f -q + q f^2+ 2fg,  \\
		v(q) &= 1 - 2 f^2 - g + qfg.
	\end{align}	
\end{corollary}

\begin{proof}
	The given descriptions for $u$ and $v$ are the unique solutions to the equations~\eqref{eq:ODE1} and~\eqref{eq:ODE2} 
	with initial conditions $u(0) = v(0) = 0$. 
\end{proof}

%% file: polynomiality.tex
\section{Polynomiality}\label{Sec:Pol}

Keeping the notations from the introduction, we want to prove the following result.

\begin{theorem}\label{Th:GenSeriesForBases}
  For any sequence of reduced partitions $\fatlambda$, 
	the generating series $$F^{odd}_{\fatlambda}(q)\quad\text{and}\quad F^{even}_{\fatlambda}(q)$$ 
	are polynomials in $q$, $f(q)$ and $g(q)$ with rational coefficients.
\end{theorem} 

Based on the invariance theorem \ref{Th:InvarTh}, we fix points $b_1<b_2<\dots<b_{m+k}$ and consider 
real simple rational functions with simple branch points $b_1,\ldots,b_m$ and 
reduced ramification profile $\lambda_1,\ldots,\lambda_k$ 
at the branch points $b_{m+1}, \dots, b_{m+k}$.
We choose an additional point $b_m < \alpha < b_{m+1}$. 

\begin{figure}[t]
\centering
\includegraphics[scale=1]{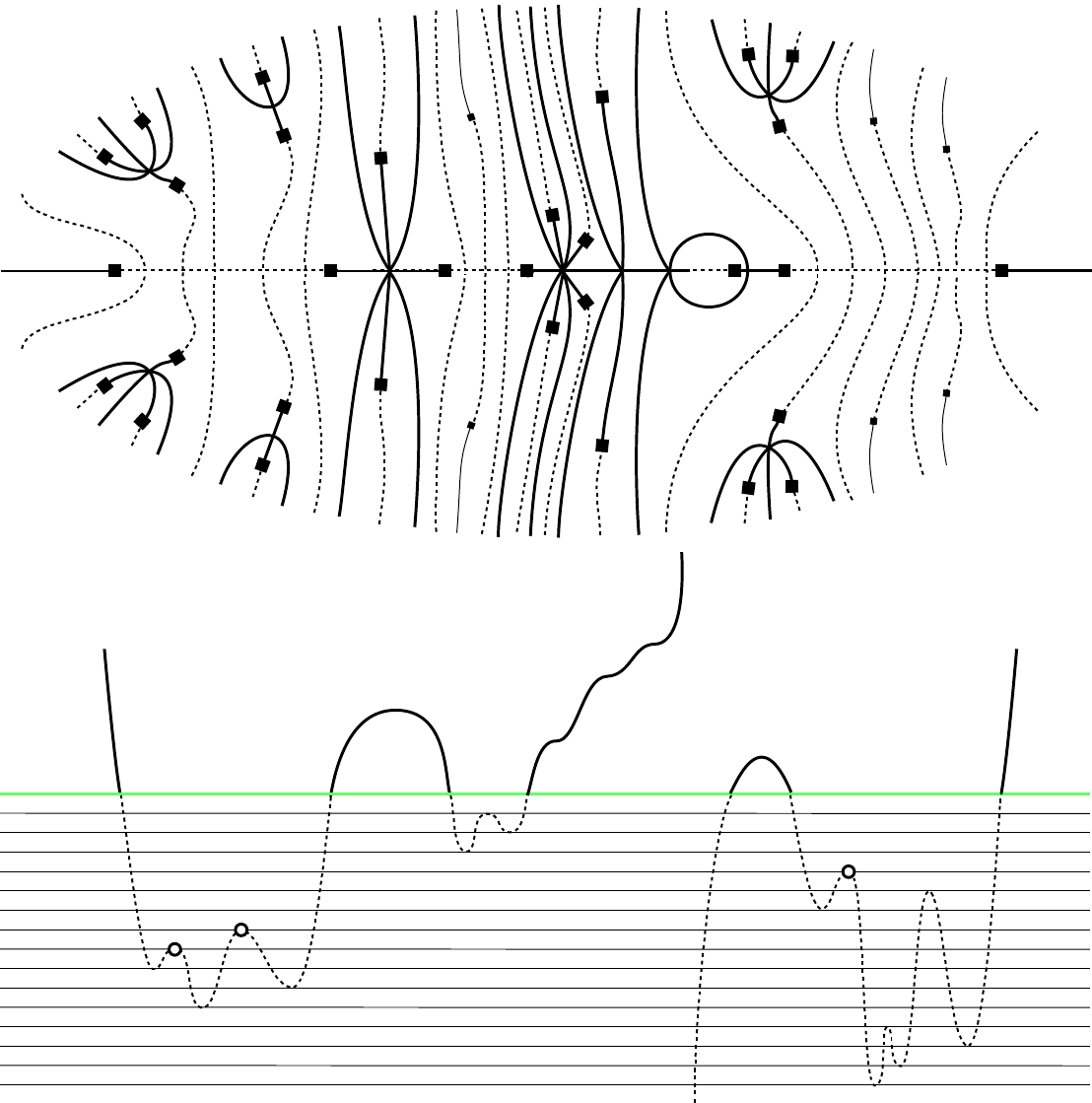}

\caption{A base $B$ of type~\ref{case1}, represented as a union of thick lines. Here, $l=4$, $\sp = 3$ and $(c_1, c_2, c_3, c_4) = (2,0,0,1)$. 
The black squares are the $\alpha$-vertices.
The dotted lines represent the inverse image of $(-\infty,\alpha)$ for a dessin $\Gamma$ with $B_\Gamma = B$.}
\label{Fig:Base(a)}
\end{figure}

Let $\varphi$ be such a function and $\Gamma$ be the associated affine dessin of $\varphi$. 
We set $B'$ to be the union of connected components of
$\varphi^{-1}([\alpha, \infty))$ which contain a critical point of $\varphi$
(see Figure \ref{Fig:Base(a)}). 
We also include the components of $\varphi^{-1}([\alpha, \infty))$
to the left and right which contain unbounded pieces of $\R$. 
Since the left one only exists if $d$ is odd, this is just a convenient way 
of storing the parity of $d$ in $B'$. 

\begin{definition} 
	A connected component $C$ of $\RR \setminus B'$
	is a \emph{chain} of $\Gamma$. 
	Each point in the inverse image of $\alpha$ is labelled by {\tiny $\blacksquare$}, 
	and is called an $\alpha$-vertex.
	We order the chains from left to right as $C_1, \dots, C_l$. 
	The \emph{special chain} $C_\sp, \sp \in \{1,\dots,l\}$
	is chain containing the simple pole $p$.
	The tuple $B_\Gamma := (B', \sp)$ is called the \emph{base} of $\Gamma$.  
\end{definition}

Let $q_1,q_2$ denote the closest critical points before and after the simple pole $p$.
We need to distinguish three cases 
(see Figures \ref{Fig:Base(a)}, \ref{Fig:SBase(b)}, and \ref{Fig:SBase(c)} 
for examples):

\begin{enumerate}[label=(\Alph*)]
	\item \label{case1} $\varphi(q_1), \varphi(q_2) > \alpha$,
	\item \label{case2} $\varphi(q_1) > \alpha > \varphi(q_2)$ or $\varphi(q_2) > \alpha > \varphi(q_1)$,
	\item \label{case3} $\varphi(q_1), \varphi(q_2) < \alpha$.
\end{enumerate}

\noindent The three cases can be distinguished in terms of the base as follows.

\begin{enumerate}[label=(\Alph*)]
	\item The base $B_\Gamma$ contains a cycle.  
	\item The base $B_\Gamma$ contains no cycle, but a bounded open end (with the simple pole $p$ as endpoint).
	\item Otherwise.
\end{enumerate}

\begin{remark}\label{Rem:Properties}
  To motivate the following definition, let us make a few simple observations.
  	Let $C_1, \dots, C_l$ be the sequence of chains, the order of which is according to the orientation of $\R$.
	Let $r_i$ denote the number of critical points occurring in $C_i$. 
	Then the parity of $r_i$ is given by  
		\begin{equation} \label{eq:parityofrnumbers1} 
		r_i \equiv_2 
		\begin{cases}
		  d & i = 1, \\
			1 & i > 1, 
		\end{cases} 		
	\end{equation}	  
	except for the following cases:
	\begin{equation} \label{eq:parityofrnumbers2}
	  r_\sp \equiv_2 
	  \begin{cases}
			d-1 	& \sp = 1, \text{and either type \ref{case1} or \ref{case2}}, \\
			0 		& \sp > 1, \text{and either type \ref{case1} or \ref{case2}}. 
	  \end{cases} 
	\end{equation}
		Moreover, we have 
	\begin{equation} \label{eq:parityofrnumbers3} 
		r_\sp = 0 \quad\quad \text { if } B \text{ is of type \ref{case1}}. 
	\end{equation}
\end{remark}

In the following, we fix a base $B=(B',\sp)$ of type \ref{case1}, \ref{case2} or \ref{case3} and chains
$C_1, \dots, C_l$. 
Note that for type \ref{case1} and \ref{case2}, 
the special index $\sp$ is uniquely determined by $B'$, while for type \ref{case3} any
$\sp \in \{1,\dots,l\}$ can occur.
Let $Z_+$ be a connected component of 
$B$ which is contained in
$\{z : \Im(z) > 0\}\subset\CC$.
It has a unique $\alpha$-vertex which can be connected in $\CC \setminus B$ 
to a unique chain $C_i$.
We say $Z_+$ is \emph{adjacent} to $C_i$ 
and denote by $c_i$ the number of such $Z_+$
(see Figure \ref{Fig:Base(a)}). 

\begin{definition}\label{Def:ChainData}
	A \emph{chain data set} for  a base $B$ consists of a tuple $(I_i, a_i, M_i)_{i \in \{1,\dots,l\}}$ 
	with the following properties.
	\begin{enumerate}
			\item 
			The sets $I_i$ form a disjoint union $I_1 \sqcup \dots \sqcup I_l = \{1, \dots, m\}$
			and	the numbers $r_i := |I_i|$ satisfy  equations
			\eqref{eq:parityofrnumbers1}, \eqref{eq:parityofrnumbers2} and \eqref{eq:parityofrnumbers3}.
		\item
			For any $i \in \{1, \dots,\hat{\sp}, \dots, l\}$, $a_i$ is an ordinary alternation
			$a_i$ of length $r_i$. Additionally, $a_\sp$ is an ordinary or broken alternation
			of length $r_\sp$ if type is \ref{case2} or \ref{case3}. For type \ref{case1}, we set
			$a_\sp = \emptyset$. 
		\item 
			We set $m_i := \left\lfloor r_i / 2 \right\rfloor$ for $i \neq \sp$.
			We set $m_\sp := \left\lfloor r_\sp / 2 \right\rfloor -1$, except for the case
			type \ref{case2}, $d$ even, and $\sp = 1$, in which we set 
			$m_\sp := \left\lfloor r_\sp / 2 \right\rfloor$. 
			Then the last ingredient is a choice of subsets $M_i \subset \{1, \dots, m_i\}$
			of size $c_i$.
	\end{enumerate}
\end{definition}

\begin{proposition}\label{Prop:ChainData}
  Fix a base $B$ as above. Then the set of dessins $\Gamma$ with $B = B_\Gamma$
	is in one-to-one correspondence with the set of chain data sets for $B$. 
\end{proposition}
\begin{proof}
We first associate a chain data set to each $\Gamma$. 
	Fix $i \in \{1, \dots, l\}$.
	We define $I_i$ such that $j \in I_i$ if and only if the simple 
	critical point with value $b_j$ is contained in $C_i$. 
  After reparametrising $I_i$ by $\{1, \dots, r_i\}$ keeping the ordering,
	the order of appearance of the critical points in $C_i$ defines a permutation
	of $\{1, \dots, r_i\}$. 
	If the type is \ref{case2}, $i = \sp$, and the simple pole is located 
	on the right hand side of $C_\sp$, we additionally compose with the permutation with
	$x \mapsto r_\sp + 1 - x$. The result is used as $a_i$. 
	Note that $m_i$ appearing in Definition~\ref{Def:ChainData} is equal
	to the number of maxima in $C_i$ for $i\neq \sp$, and equal to the number of maxima minus
	one if $i = \sp$ and the type is \ref{case2} or \ref{case3}.
	In this case, there is a unique 
	maximum in $C_\sp$ which is closest to the simple pole $p$
	(see Figure \ref{Fig:GenerTildeG}).
	After removing this maximum, the remaining maxima can be ordered increasingly and
	labelled by $\{1, \dots, m_i\}$. 
  Let $Z_+$ be a connected component of $B$, as before, adjacent to $C_i$. Since
	$\Gamma$ is connected, there is path connecting $Z_+$ and $\RR$. 
	It is easy to see that this path is unique, its starting point is
	a $\alpha$-vertex of $Z_+$ and its endpoint is a maximum of $C_i$. 
  Varying $Z_+$, this construction distinguishes $c_i$ maxima of $C_i$,
	whose labels we collect in $M_i \subset \{1, \dots, m_i\}$.
	By the properties of $\Gamma$ it is clear that we constructed 
	a well-defined chain data set for $B$.

\begin{figure}[h]
\centering
\includegraphics[scale=1]{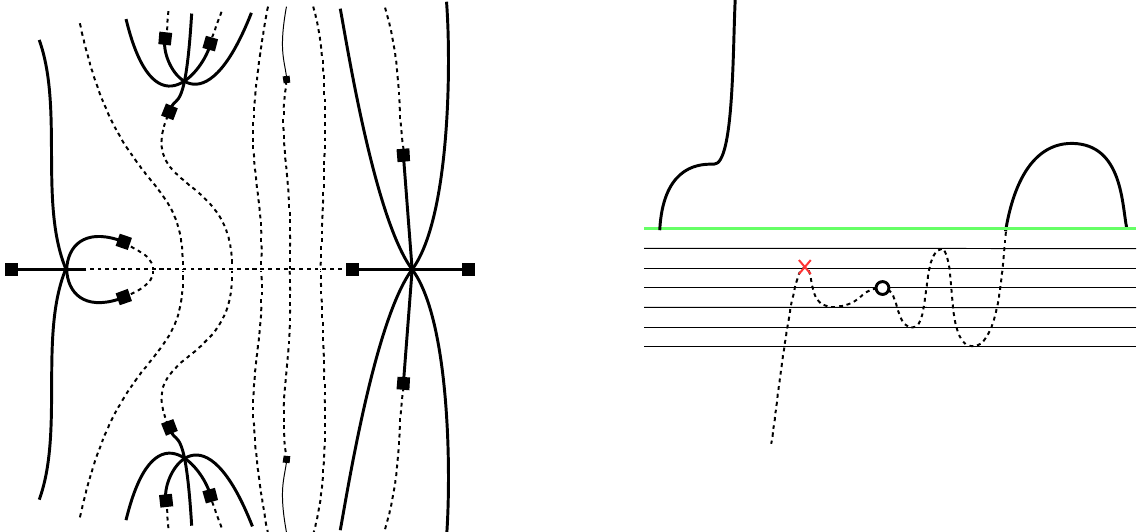}
\caption{The special chain of a base of type \ref{case2}. The ticked maximum cannot be marked in $M_\sp$.}
\label{Fig:GenerTildeG}
\end{figure}

	We proceed to show the existence of the inverse map.
	Assume we are given a chain data set for $B$.
	Fix some $i \in \{1, \dots, l\}$. 
	Using the bijections in Proposition~\ref{Prop:AlterDessins} and~\cite[Example 1.5]{IZ18},
	we can associate to $a_i$ a unique (affine) real dessin $D_i$ with labels $1 \to \dots \to r_i$.
	Reversing the relabelling from above, we can rename the labels by elements of $I_i$.
	We complete this to a labelling of type $-\infty \to 1 \to \dots \to m \to +\infty$ by
	adding two-valent vertices as necessary. 
	The next step is to glue the $D_i$ into $\CC$ such that their real parts go to 
	$C_i$ keeping the orientation. 
	We do this such that $\RR$ is fully covered by $B$ and the $D_i$ (in type \ref{case1} and \ref{case2},
	we have to throw in by hand the point $p$ which closes up the open bounded end of $B$).
	The next step is to connect some non-real ends of $D_i$ to $B$.
	Of course, it suffices to describe the construction in the upper half plane,
	and then copy it symmetrically. 
	We start with type \ref{case2} and $i = \sp$. 
	In this case there is a unique maximum in $D_\sp$ closest to $p$,
	and an end of $D_\sp$ emanating from it. There is a unique $\alpha$-vertex
	contained in the connected component of $B$ adjacent to $p$ and to which 
	this end of $D_i$ can be glued. We connect them.
	Among the $m_i$ maxima of $D_i$ which are not ``adjacent'' to $p$,
	we select the subset of size $c_i$ corresponding to $M_i$. 
	We consider the $c_i$ ends of $D_i$ which emanate 
	from the these maxima. There is a unique way (up to homeomorphism) to glue 
	these ends to $\alpha$-vertices of the connected components of $B$
	adjacent to $C_i$ without producing extra intersections. Again, we connect them like this.
	We now complete the dessins as follows. Each end of $D_i$ of type $m \to \infty$
	not yet glued to a $\alpha$-vertex is replaced by a sequence of edges
	$m \to m+1 \to \dots m+k \to +\infty$. This includes, in type \ref{case3}, 
	the end adjacent to the simple pole in $D_\sp$. 
	Finally, all non-real ends $-\infty \to 1$ and $m+k \to +\infty$ are
	extended to infinity to obtain an honest affine dessin $\Gamma$. 
	It is clear that this is indeed the inverse map we are looking for.
\end{proof}

Using this Proposition, we can explicitly describe the generating series of $S$-numbers
for a given base $B$. Similar to dessins, we consider two bases isomorphic, denoted by $B_1 \cong B_2$,
if they can be identified under a $\conj$-equivariant homeomorphism of $\CC$ that 
is orientation-preserving on $\RR$. 

\begin{definition}\label{Not:GenSeriesFB}
Given a base $B=(B',\sp)$, 
we denote by $\SSS_B(\PPP, \Fatlambda)$ the set of real simple rational dessins $\Gamma$
with $B_\Gamma \cong B$.
We set $S_B(\fatlambda,m) := S(\SSS_B(\PPP, \Fatlambda))$ and 
\[
  F_B:=\sum_{m\geq 0}S_B(\fatlambda,m)\frac{q^m}{m!}.
\]
\end{definition}

\begin{remark}\label{Rem:SumOverBases}
  Clearly, after fixing the parity of $d$ and $\fatlambda$, 
	there is only a finite collection of isomorphism classes of bases that occur
  as $B_\Gamma$. Hence we can express $F^\odd_\fatlambda$ and $F^\even_\fatlambda$ as finite sums
	\[
	  F^\odd_\fatlambda = \sum_B F_B \quad \text{ and } \quad F^\even_\fatlambda = \sum_B F_B.
	\]
\end{remark}

\begin{definition} \label{def_signbase}
  Let $B$ be a base and let $\Sigma_j$ be the sequence of degrees of
	real vertices of $B$ contributing to $\lambda_j$ (i.e., the vertices
	with label $m+j$), with an extra $1$ inserted at the position of $C_\sp$
	\emph{only} if $B$ is of type \ref{case1} or \ref{case2}.
	The \emph{sign} of a base $B$ is 
  \[
		\varepsilon(B) := (-1)^{\dis \Sigma_1+ \dots + \dis \Sigma_k}.
	\]
\end{definition}

\begin{remark} \label{rem_signcaseC}
	In case \ref{case3}, two consecutive $1$'s should be inserted, one for the simple pole $p$,
	and one for the $2$-valent vertex on the segment in $\varphi^{-1}([\alpha, \infty))$
	adjacent to $p$, which is not part of $B$. 
	In the definition, equivalently, we do not insert anything. 
\end{remark}

We denote the differential operator $D=q\frac{d}{dq}$. 
We introduce the following functions 
  \begin{align*} 
		f_c 				&= \frac{1}{2^c c!}(D-1)(D-3)\cdots (D-2c+1)f, \\
		g_c 				&= \frac{1}{2^c c!}D(D-2)\cdots (D-2c+2)g, \\
		\tilde{g}_c	&= \frac{1}{2^c c!}(D-2)(D-4)\cdots (D-2c)g, \\
		u_c					&= \frac{1}{2^c c!}(D-3)(D-5)\cdots (D-2c-1)u, \\
		v_c					&= \frac{1}{2^c c!}(D-2)(D-4)\cdots (D-2c)v.
	\end{align*}

\begin{theorem}\label{Th:ExplicitExpress} 
Fix a base $B$ and set $\varepsilon = \varepsilon(B)$.
Then, depending on the different cases, the generating series $F_B$ can be written as in the following table.
	\[
		\begin{array}{c|c|c|c|c}
			F_B &   \text{type \ref{case1}} &  & \text{type \ref{case2}} & \text{type \ref{case3}} \\ \cline{1-2}\cline{4-5}
			d \text{ odd} & \displaystyle \varepsilon \, \prod_{\substack{i=1 \\ i\neq\sp}}^l f_{\displaystyle c_i} & 
										& \displaystyle \varepsilon \, \tilde{g}_{\displaystyle c_{\sp}} \, \prod_{\substack{i=1 \\ i\neq\sp}}^l f_{\displaystyle c_i} 
										& \displaystyle \varepsilon \, u_{\displaystyle c_s} \, \prod_{\substack{i=1 \\ i\neq\sp}}^{l} f_{\displaystyle c_i}\\[5pt]
			\hline
			\multirow{4}{*}{d \text{ even}} & \multirow{4}{*}{ $\displaystyle \varepsilon \, g_{\displaystyle c_1} \, \prod_{\substack{i=1 \\ i\neq\sp}}^{l}f_{\displaystyle c_i} $} 
						& \sp =1   & \displaystyle - \varepsilon \, f_{\displaystyle c_1} \, \prod_{i=2}^{l} f_{\displaystyle c_i}  
										& \displaystyle \varepsilon \, v_{\displaystyle c_1} \, \prod_{i=2}^{l} f_{\displaystyle c_i} \\\cline{3-5} 
						& & \sp \neq 1 & \displaystyle \varepsilon \, g_{\displaystyle c_1} \, \tilde{g}_{\displaystyle c_{\sp}} \prod_{\substack{i=2 \\ i\neq\sp}}^{l} f_{\displaystyle c_i} 
										& \displaystyle \varepsilon \, g_{\displaystyle c_1} \, u_{\displaystyle c_\sp} \, \prod_{\substack{i=2 \\ i\neq\sp}}^{l}f_{\displaystyle c_i} \\ \hline
		\end{array}
	\]
\end{theorem}

\begin{proof}
Up to signs, the statement follows directly from Proposition \ref{Prop:ChainData}. 
Indeed, let us compare to the ingredients of a chain data set from Definition \ref{Def:ChainData}.
Ingredient (b) corresponds to the 
coefficients in $f, g, u, v$ by Proposition \ref{Prop:AlterDessins} and Remark 
\ref{rem_ordinaryalternations},
ingredient (c) corresponds to the extra factors in the coefficients of 
$f_c, g_c, \tilde{g}_c, u_c, v_c$, compared to the non-indexed counterparts,
and finally ingredient (a) corresponds to taking the product of the generating series.
For the sign computations, first note that all disorders $\dis \Sigma_{b_j}, j = m+1, \dots, m+k$ are accounted for in $\varepsilon(B)$
(see also Remark \ref{rem_signcaseC}). 
For any simple branch point $b_j, j = 1, \dots, m$, the parity of $\dis \Sigma_{b_j}$ only depends on whether the corresponding 
critical point is a local maximum or minimum, which is compatible with the signs in $f,g,u,v$ by 
Proposition \ref{Prop:AlterDessins}. There is one exception, namely if the special chain for type \ref{case2}
is of the form $[a,b)$ (pole on the right), in which case an extra reflection was necessary.
There are two subcases: For $\sp \neq 1$ the number of maxima and minima in $D_\sp$ is equal, 
hence the the sign is not affected by the reflection. For $\sp = 1$, however, the number of maxima and minima
differs by one, which is fixed by the extra minus sign in this case.
\end{proof}

\begin{remark}\label{L:OperatorPolynom}
Recall from \cite[Lemma 5.3]{IZ18} that 
\begin{equation} \label{eq_derivativesqfg} 
    Dq=q,\quad Df=q(1-f^2),\quad\text{and}\quad Dg=-qfg.
\end{equation}
In particular, if $h\in\QQ[q,f]$, then $Dh \in \QQ[q,f]$, 
and if $h\in\QQ[q,f]g$, then $Dh \in \QQ[q,f]g$.
\end{remark}

\begin{proof}[Proof of Theorem~\ref{Th:GenSeriesForBases}]
By Corollary \ref{uvpolynomialinfg} and Remark \ref{L:OperatorPolynom}, 
we have $f_c, g_c, \tilde{g}_c, u_c, v_c \in \QQ[q,f,g]$. 
Moreover, by Theorem~\ref{Th:ExplicitExpress}, we have that $F_B \in \QQ[q,f,g]$. 
Therefore, the result follows from Remark~\ref{Rem:SumOverBases}.
\end{proof}

%% file: nonvanishing.tex
\section{Non-vanishing statements for generating series}\label{Sec:Non-Van}

The functions $f$ and $g$ satisfy the relation $f^2 + g^2 = 1$, 
which implies
\[
  \QQ[q,f,g] = \QQ[q,f] + \QQ[q,f]g.
\] 
The corresponding decomposition of an element in $\QQ[q,f,g]$,
removing higher powers of $g$, is called \emph{$g$-reduced}.
The following lemma asserts uniqueness of this representation. 

\begin{proposition}\label{Prop:NotAllZero}
  The family of power series $q^i f^j$, $q^i f^j g$, $i,j \in \NN \cup \{0\}$ is linearly independent 
	over $\CC$.
\end{proposition}

\begin{proof}
  Regarded as meromorphic functions, $f=\tanh$ and $g=\sech$ have the same set of poles
	$p_k := (1/2 + k) i \pi, k \in \ZZ$, all of them simple. Moreover, a simple calculation shows
	$\Res_{p_k}(f) = i$ and $\Res_{p_k}(g) = (-1)^{k-1} i$.
	We consider an equality of power series or, equivalently, meromorphic functions
	\[
	  0 = a + \sum_{j=1}^\delta (a_j f + b_j g)f^{j-1},
	\]
  with $a,a_j,b_j \in \CC[q]$. 
	Looking at the highest order term, we can observe that all poles in $a_\delta f + b_\delta g$
	need to cancel out, since otherwise the expression on the right hand side has a pole of order $\delta$.
	But this implies
	\[
	  0 = \Res_{p_k}(a_\delta f + b_\delta g) = i (a_\delta + (-1)^{k-1} b_\delta)(p_k)
	\]
	for all $k \in \ZZ$. This implies that both polynomials $a_\delta + b_\delta$ and $a_\delta - b_\delta$ have
	infinitely many zeros, and hence $a_\delta = b_\delta = 0$.
	Proceeding recursively, we deduce that all the polynomials $a_j,b_j$ and finally $a$ are zero.
	\end{proof} 

From now on, we make constant use of the proposition and always replace 
elements in $\QQ[q,f,g]$ by their unique $g$-reduced representation. 

Set $\deg(q^i f^j) = (i,j)$ and extend to a degree function on $\QQ[q,f]$ using the lexicographic order, 
i.e.\ $(i,j) < (i',j')$ if and only if $i < i'$ or $i=i'$ and $j < j'$. 
For an element $F = F_1 + F_2 g \in \QQ[q,f] \oplus \QQ[q,f]g$, we define the \emph{$f$-degree} $\deg_f(F) := \deg(F_1)$ and
the \emph{$g$-degree} $\deg_g(F) := \deg(F_2)+(0,1)$ if $F_2\neq 0$, otherwise $\deg_g(F)=0$. 
Note that for $F,F' \in \QQ[q,f]g$, we have $\deg_f(F \cdot F') = \deg_g(F) + \deg_g(F')$.
The derivation rules in Equation \eqref{eq_derivativesqfg} imply the following statement.

\begin{corollary} \label{cor_degreederivative}
  For any $F \in \QQ[q,f]f \oplus \QQ[q,f]g$ we have
	\[
	  \deg_f(D F) = \deg_f(F) + (1,1), \quad \quad  \deg_g(D F) = \deg_g(F) + (1,1).
	\]
  Moreover, the leading coefficient of $DF$ with respect to $f$-degree and $g$-degree 
	is $-j$	times the corresponding leading coefficient of $F$, 
	where $\deg_f(F)=(i,j)$, respectively, $\deg_g(F)=(i,j)$.	
	\end{corollary}

\begin{proposition}\label{Prop:fAndgDegrees}
  Let $B$ be a base with $l$ chains and let $c = \sum_{i=1}^l c_i$ be the total number of connected
	components of $B$ contained in $\{\Im(z)>0\}$.
	Then $\deg(F_B)$ takes value as described in the following table.
	\[
		\begin{array}{cccll}
			\hline
										&             & \text{type \ref{case1}} & \text{type \ref{case2}} & \text{type \ref{case3}} \\ \hline
			d \text{ odd} & \deg_f(F_B) & (c,c+l-1)               & (0,0)                   & (c+1,c+l+1)             \\
										& \deg_g(F_B) & (0,0)                   & (c,c+l)                 & (c,c+l+1) \; \; \; \text{\textbf{[1]}}           \\ \hline
			d \text{ even}& \deg_f(F_B) & (0,0)                   & (c,c+l) \; \; \; \text{\textbf{[2]}}  & (c,c+l+1) \; \; \; \text{\textbf{[3]}}  \\
										& \deg_g(F_B) & (c,c+l-1)               & (0,0)                   & (c+1,c+l+1)             \\ \hline
		\end{array}
	\]
	Moreover, the coefficients of the leading terms are $\varepsilon(B) (-1)^c / 2^c$, except 
	for the cases marked by \textbf{[1]}, \textbf{[2]}, \textbf{[3]}, which have leading coefficients
	\[
	  \text{\textbf{[1]}} \;\; \varepsilon(B) (-1)^c / 2^{c-1} \quad\quad
	  \text{\textbf{[2]}} \;\; \varepsilon(B) (-1)^{c-1} / 2^c \quad\quad
	  \text{\textbf{[3]}} \;\; \varepsilon(B) (-1)^{c-1} / 2^{c-1}.
	\]
\end{proposition}
	
\begin{proof}
  Applying Corollary \ref{cor_degreederivative}, we can compute the degrees of the modified generating functions
	as described in the following table.
	\[
		\begin{array}{cccccc}
			\hline
							& f_n & g_n & \tilde{g}_n & u_n & v_n \\
			\hline
			\deg_f 	& (n,n+1) & (0,0) & (0,0) & (n+1,n+2) & (n,n+2) \\
			\deg_g 	& (0,0) & (n,n+1) & (n,n+1) & (n,n+2) & (n+1,n+2) \\
			\hline
		\end{array}
	\]	
  Then the statement about degrees follows from Theorem \ref{Th:ExplicitExpress}.
	Similarly, by Corollary \ref{cor_degreederivative} the various leading coefficients of these 
	functions are $(-1)^n/2^n$, except for the following special cases.
	\begin{enumerate}
		\item[\textbf{[1]}] The term $u_{c_\sp}$ carries an extra factor of $2$ by Corollary \ref{uvpolynomialinfg}.
		\item[\textbf{[2]}] The terms $-f_{c_1}$ and $g_{c_1} \tilde{g}_{c_\sp}$ produce an extra $-1$.
		\item[\textbf{[3]}] The terms $v_{c_1}$ and $g_{c_1} u_{c_\sp}$ carry an extra factor of $-2$ by Corollary \ref{uvpolynomialinfg}.
	\end{enumerate}
	Hence the statement about leading coefficients follows.
	\end{proof}

\begin{notation}\label{Def:NumbersInLevels}
For a fixed set of reduced partitions $\fatlambda:=(\lambda_1, \dots, \lambda_k)$, denote by $\cc$ the number of pairs in each partition of $\fatlambda$, and by $\oo$ (resp. $\ee$) the number of partitions in $\fatlambda$ with an odd (resp.\ even) entry appearing an odd number of times.
\end{notation}

From now on, we fix the parity of $d$ and a sequence of reduced ramification profiles $\fatlambda$
such that the non-vanishing criteria in Theorem \ref{Th:IffOdd} are satisfied.
We denote by $\oo$ (resp. $\ee$) the number of partitions in $\fatlambda$ 
with an odd (resp.\ even) entry appearing an odd number of times.
We denote by $\cc$ the number of pairs of equal entries in each partition of $\fatlambda$, 
such that $2 \cc + \oo + \ee = \sum_j l(\lambda_j)$.

\begin{theorem}\label{Cor:DegreeLessOrEqual}	
  The degrees of $F^\odd_\fatlambda$ and $F^\even_\fatlambda$ can be bounded by
	\begin{align} 
		\deg_f F^\odd_\fatlambda  &\leq (\cc+1,\cc+\oo+2),   &\deg_g F^\odd_\fatlambda  &\leq (\cc,\cc+\oo+2), \\
		\deg_f F^\even_\fatlambda  &\leq (\cc,\cc+\oo+2),    &\deg_g F^\even_\fatlambda  &\leq (\cc+1,\cc+\oo+2).
	\end{align}
\end{theorem}

\begin{proof}
  Let $B$ be a base of given type. We note that $c \leq \cc$, since any connected component contained in
	$\{\Im(z)>0\}$ accounts for at least one pair of complex conjugated critical points.
	Assume that $B$ is a base with $c = \cc$. Any connected component $B \cap \RR$ of the form $[a,b]$
	(a bounded closed interval) contains at least one vertex of degree $4n$, hence the number of such components is bounded
	by $\oo$. It follows that the number of chains $l$ is bounded by $l \leq \oo + 2$ for type \ref{case1} and \ref{case2}
	and by $l \leq \oo + 1 $ for type \ref{case3}. Then the claim follows from Proposition \ref{Prop:fAndgDegrees}.
\end{proof}

Our goal in the remaining subsections is to prove sharpness of these estimates 
(in same cases) by proving that the corresponding leading coefficient is non-zero.

\begin{definition}\label{Def:SimpleBases}
A \emph{simple base} $B$ is a base of type \ref{case2} or \ref{case3} such that
		\begin{itemize}
		\item the number of connected components contained in $\{\Im(z)>0\}$ is $\cc$,
		\item the number of bounded closed connected components in $B \cap \RR$ is $\oo$,
	\end{itemize}
\end{definition}

\begin{figure}[b]
\centering
\includegraphics[scale=1]{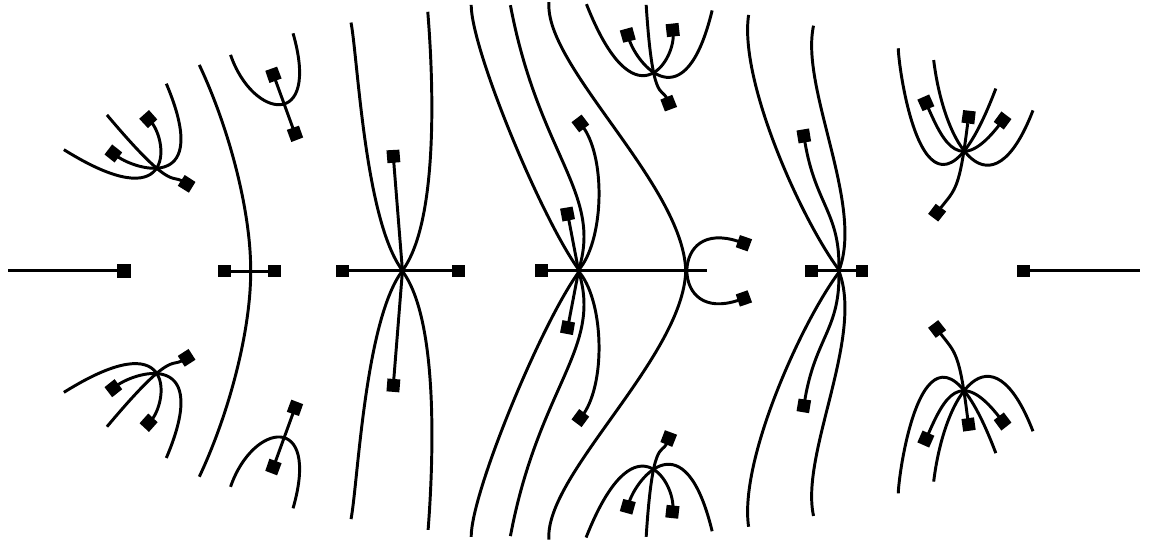} 
\caption{A simple base of type \ref{case2} with $\cc = 4$, $\oo=3$ and $\ee = 2$.}
\label{Fig:SBase(b)}
\end{figure}

\begin{figure}
\centering
\includegraphics[scale=1]{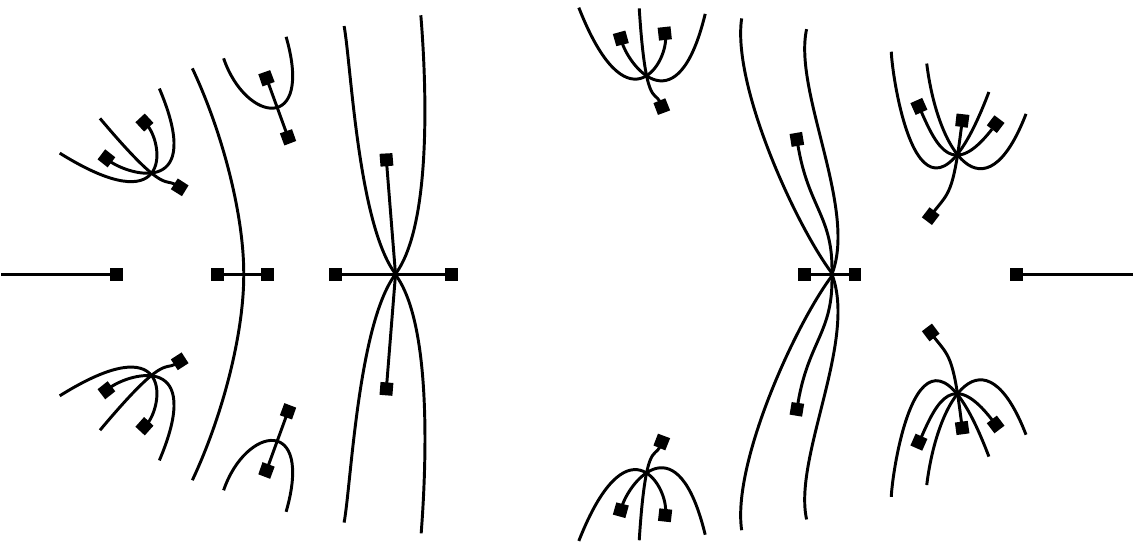}
\caption{A simple base of type \ref{case3} with $\cc = 4$, $\oo=3$ and $\ee = 0$.}
\label{Fig:SBase(c)}
\end{figure}

\begin{remark}\label{Rem:SimpleTypeB}
  Simple bases are exactly the maximal bases appearing in the proof of Theorem \ref{Cor:DegreeLessOrEqual}. 
	Each connected component $Z_+$ contained in $\{\Im(z)>0\}$ contains exactly one vertex of degree greater than $2$.
	Each bounded closed connected component of $B \cap \RR$ is $\oo$ contains exactly one vertex of degree
	$4n$ and no two of these vertices carry the same label. They correspond to local maxima and are called the 
	\emph{maxima} of $B$. 
	There are $\ee$ vertices of degree $4n + 2, n > 0$, which we call the \emph{crossings}
	of $B$. Again, the labels of crossings are pairwise distinct.
	In type \ref{case2}, there exists exactly one finite half-closed connected component of $B \cap \RR$, 
	and it contains no maximum, but at least one crossing. In particular, if $\ee = 0$, 
	simple bases of type~\ref{case2} do not exist. 
  However, simple bases of type \ref{case2} do exist if $\ee > 0$, 
	as well as simple bases of type \ref{case3} for any value of $\ee$
	(cf.\ their polynomial counterparts defined in~\cite[Section 5]{IZ18}).
	Examples are given in Figures \ref{Fig:SBase(b)} and \ref{Fig:SBase(c)}.
	\end{remark}

It is now easy to the describe the coefficients in front of the monomials corresponding to
the upper bounds in Proposition~\ref{Prop:fAndgDegrees}. 
Let $\SBB$ and $\SBC$ denote the set of simple bases of type \ref{case2}
and type \ref{case3}, respectively. We defined signs on these sets in Definition
\ref{def_signbase} and can take the signed counts $S(\SBB)$ and $S(\SBC)$. 

\begin{corollary}\label{Th:Coeffs}
  The coefficients of the monomials corresponding to
	$f$-degree $(\cc+1,\cc+\oo+2)$, $g$-degree $(\cc,\cc+\oo+2)$, for $d$ odd,
	and $f$-degree $(\cc,\cc+\oo+2)$, $g$-degree $(\cc+1,\cc+\oo+2)$, for $d$ even, are
	equal, in the same order, to $\frac{(-1)^\cc}{2^\cc}$ times
	
	\[
	  S(\SBC), \quad S(\SBB) + 2 S(\SBC), \quad -S(\SBB) - 2 S(\SBC), \quad S(\SBC).
	\]
\end{corollary}

\begin{proof}
  The statement follows by Proposition \ref{Prop:fAndgDegrees}, Theorem \ref{Cor:DegreeLessOrEqual}, and Remark \ref{Rem:SimpleTypeB}.
\end{proof}

\subsection{The case \texorpdfstring{$\ee = 0$}{e = 0}} 

\begin{lemma}\label{L:NoCrossings}
  If $\ee = 0$, no simple bases of type \ref{case2} exist, and all simple bases 
	of type \ref{case3} have the same sign $(-1)^\oo$. 
\end{lemma}

\begin{proof}
The non-existence of bases of type \ref{case2} is explained in Remark \ref{Rem:SimpleTypeB}.
Assume $B$ is a simple base of type \ref{case3}. Since for each label $j$ there exists at most
one real vertex with label $j$, and all of them are maxima (i.e., the non-one entries in $\Sigma_{b_j}$ are succeeded by
an odd number of $1$'s), each maximum contributes $-1$ to $\varepsilon(B)$, and the claim follows.
\end{proof}

\subsection{The case \texorpdfstring{$\ee > 0$}{e > 0}}

\begin{definition}\label{Def:BasesEquiv}
Given a simple bases $B$, a crossings $v$ and a $2$-valent vertex $w$ of the same label, we denote by $\Flip^{v}_{w}(B)$ the base 
obtained from $B$ by exchanging the two vertices together with adjacent edges of $v$ in $\CC \setminus \RR$. 
If we flip an \enquote{increasing} and a \enquote{decreasing} vertex, we also reverse the order of these edges.  
We denote by $[B]$ the equivalence class of simple bases up 
to changing the position of the crossings via the flip operation.
\end{definition}

\begin{proposition}\label{Prop:SignsOfBases}
If $\ee$ is a positive integer, then
	\[  
	  S([B]) = \begin{cases}
		  0 								& d \text{ odd, type \ref{case3} ,} \\
		  (-1)^{\oo + \bb} 	& d \text{ odd, type \ref{case2}, } \ee \text { even, or } d \text{ even, type \ref{case3}.} 
		\end{cases}
	\]
	Here, $\bb$ denotes the number of partitions $\lambda_i$ with both an odd and even entry appearing an odd number of times, 
	and the even entry is bigger than the odd entry.
\end{proposition}

Note that we neglect the case $d$ even, type \ref{case2}, which is not needed here.

\begin{proof}
  The connected components of $B \cap \RR$ of form $(-\infty, y], [x, +\infty), [x,y]$ and $[x,y)$ or $(x,y]$ 
	are called left end, right end, ordinary segment, special segment in the following. 
	For each crossing, there are zero or two possible positions on an ordinary segment and
	one position on any other segment. Let $v$ be a crossing (labelled by $j$) sitting on the left end (only if $d$ is odd)
	or on an ordinary segment. Let $w$ be the $2$-valent vertex of the same label on the right end or 
	on the same ordinary segment. Then $\Flip^{v}_{w}(B)$ has opposite sign. This is clear in the second case
	and easy to verify in the first case, since the sequence $\Sigma_j$ we are manipulating is of the form
	\[
	  (a, \underbrace{1, \dots, 1}_{\text{even number}}, 1) \quad \text{ or } \quad (a, \underbrace{1, \dots,1, b, 1 \dots 1}_{\text{odd number}}, 1),
	\]
	with $1 \neq a \neq b \neq 1$. The second case occurs if and only if there is a maximum with same label as $v$.  
	We can easily turn this into an involution on $[B]$, reducing the computation of $S([B])$ to those bases
	for which all crossings are located on the special segment, if $B$ is of type \ref{case2}, and 
	all crossings located on the right end, if $B$ is of type \ref{case3} and $d$ is even. \\
	It remains to show that the sign of such bases is $(-1)^{\oo + \bb}$. 
	In the first case (type \ref{case2}, crossings on special segment), 
	we note that a crossing labelled with $j$ is succeeded by 
  an odd or even number of $1$'s in $\Sigma_j$, depending on whether the special segment
	is of the form $(x,y]$ or $[x,y)$. Since $\ee$ is even, the total contribution of
	the corresponding disorders is always even. 
	It remains to count the disorders involving a maximum.
	A maximum labelled by $j$ is succeeded by an odd number $l_j$ of entries in $\Sigma_j$.
	If no crossing of label $j$ and higher degree than the maximum exists, the number 
	of disorders involving the maximum is $l_j$. 
	If a crossing of higher degree than the maximum exists, then the number of disorders
	involving the maximum is $l_j + 1$ or $l_j -1$, depending on whether
	the crossing lies before or after the maximum. 
	Hence the critical levels to be counted 
	are exactly those which contain a maximum, but no crossing of higher degree. 
	The number of such levels is $\oo - \bb$.
\end{proof}

\subsection{Proofs of main theorems}

\begin{proof}[Proof of Theorem~\ref{thm_nonvanishing}]
  The vanishing statements are contained in Propositions \ref{Th:OnlyIfEven} and \ref{Prop:OnlyIfOdd}, and it remains
	to prove non-vanishing under the given conditions.
  Fix $\fatlambda$ and the parity of $d$
	such that the conditions in Theorem \ref{thm_nonvanishing} are satisfied.
	Recall that Theorem~\ref{Th:GenSeriesForBases} asserts that the series 
	$F^{odd}_{\fatlambda}(q)$ and $F^{even}_{\fatlambda}(q)$ are polynomials in $q$, $f$ and $g$.
	Therefore, by Proposition~\ref{Prop:NotAllZero}, it suffices to prove that 
	in the $g$-reduced representations of $F^{odd}_{\fatlambda}(q)$ and $F^{even}_{\fatlambda}(q)$
	at least one of the coefficients is non-zero. 
	
	\emph{$d$ odd, $\ee = 0$:} 
	By Corollary \ref{Th:Coeffs} it suffices to show $S(\SBC) \neq 0$.
	This follows by 
	Lemma \ref{L:NoCrossings}, 
	since all elements in $\SBC$ have the same sign. 
	
	\emph{$d$ odd, $\ee > 0$:} 
	By Corollary \ref{Th:Coeffs} it suffices to show $S(\SBB) + 2 S(\SBC) \neq 0$.
	This follows by 
	Proposition \ref{Prop:SignsOfBases}, 
	since $S(\SBC) = 0$ and $S(\SBB) = (-1)^{\oo + \bb} E$, 
	where $E$ denotes the number of equivalence classes in $\SBB$. 
	
	\emph{$d$ even:} 
	By Corollary \ref{Th:Coeffs} it suffices 
	to show $S(\SBC) \neq 0$.
	This follows by Propositions \ref{L:NoCrossings} and \ref{Prop:SignsOfBases}, since $S(\SBC) = (-1)^{\oo + \bb} E$, 
	where $E$ denotes the number of equivalence classes in $\SBC$. 
	This completes the proof. 
\end{proof}

\begin{remark}\label{Rem:Coeffs}
To compute the numbers $S(\SBB)$ and $S(\SBC)$ appearing in the previous proof, we note that simple bases belonging 
to different equivalence classes differ by the position of $\sp$, and the ordering of the pairs of 
complex conjugate critical points and real critical maxima. Moreover, we need to avoid overcounting 
bases having two critical points with same multiplicity and label in $\{\Im(z)>0\}$. 
Denote by $N_j^i$ the number of times an element $i$ appears in $\lambda_j$ and
set 
\[
  A(\fatlambda) = \prod_{j=1}^k \prod_i \lfloor N_j^i/2\rfloor,
\]
where the second product runs through all $i$ with $N_j^i \geq 2$. 
Then for all $d$ we have 
$$\SBC=\frac{\displaystyle(\oo+ 1) \cdot (\oo+\cc)! }{A(\fatlambda)},$$ 
and for $d$ odd and $\ee>0$ we have 
$$\SBB=\frac{\displaystyle 2\cdot(\oo+\cc+1)!}{A(\fatlambda)}.$$
\end{remark}

\begin{proof}[Proof of Theorem~\ref{thm_LogGrowth}]
  Since we assume that the non-vanishing criteria from Theorem \ref{thm_nonvanishing} are satisfied,
	we have $F:=F^{odd}_{\fatlambda} \neq 0$ or $F:=F^{even}_{\fatlambda} \neq 0$, respectively. 
	By the proof of Proposition~\ref{Prop:NotAllZero} we conclude that $F$ 
	has bounded convergence radius $\pi/2$, with poles $\pm i\pi/2$
	on the boundary circle. Note that $F$ is either an even or an odd function, depending on the parity of $\sum_{i,j} \lambda_j^i$.
	Following \cite[Proof of Theorem 5]{IZ18}, 
	let $G$ denote the function obtained by dividing by $q$, if $F$ is odd, and performing 
	the variable change $Q = q^2$. Then $G$ has a unique singularity at radius $\pi^2/4$ and
	hence $\ln |b_k| \sim - k \ln(\pi^2/4)$, where $b_k$ denote the coefficients in $G$. 
	Since $b_k = S(\fatlambda, m) / m!$ for $m = 2k$ or $m = 2k +1$, respectively, the claim follows.
\end{proof}

\begin{remark} \label{rem_complexcount}
  Let $H^\CC(\fatlambda, m)$ denote the Hurwitz number counting
	complex simple rational functions $\varphi \in \CC(z)$
	with $k$ critical levels of reduced ramification type $\lambda_1, \dots, \lambda_k$ and 
	$m$ additional simple branch points.
	We would like to show $\ln H^\CC(\fatlambda, m) \sim m \ln(m)$ for $m \to \infty$, 
	and hence	$\ln H^\CC(\fatlambda, m) \sim \ln |S(\fatlambda, m)|$ under the non-vanishing assumption
	of Theorem~\ref{thm_LogGrowth}.
	
	To prove the claim, we first use the perturbation argument from \cite[Proof of Theorem 5.10]{Rau:LowerBounds}
	to show that 
	\[
	  H^\CC(\fatlambda, m) \leq H^\CC(\emptyset, m + \sum_{i,j} \lambda_j^i).
	\]
	The right hand side can be computed using the classical formula by Hurwitz for genus $0$ single Hurwitz numbers,
	see \cite[page 22]{Hurwitz}, which gives for $m > 1$
	\[
		H^\CC(\emptyset, m) = (m-1)^{m-1}.
	\]
	Therefore the asymptotics of $\ln H^\CC(\fatlambda, m)$ is bounded from above by 
	\[
	  \ln H^\CC(\emptyset, m) \sim m \ln(m).
	\]
	Then the equivalence follows from a suitable lower bound, e.g.\ the real count, 
	\[
	  |S(\fatlambda, m)| \leq 2 H^\CC(\fatlambda, m).
	\]
	The factor $2$ is due to the fact that in the real case we consider oriented functions. 
\end{remark}